\documentclass{IEEEtran}

\usepackage{enumitem,kantlipsum}
\usepackage{graphicx}
\usepackage{epstopdf}
\usepackage{amssymb}
\usepackage{amsmath}
\usepackage{subfigure}
\usepackage{epsfig}
\usepackage{color}
\usepackage{enumitem}
\usepackage{float}
\usepackage{soul}
\usepackage{wrapfig}
\usepackage{arydshln}
\usepackage{tablefootnote}

\usepackage{amsthm}

\usepackage[numbers,sort&compress]{natbib}
\def\0{{\bf 0}}
\def\1{{\bf 1}}

\def\beq{\begin{equation*}}
    \def\eeq{\end{equation*}}
\def\bql{\begin{equation}}
    \def\eql{\end{equation}}
\def\bqn{\begin{eqnarray*}}
    \def\eqn{\end{eqnarray*}}
\def\bnl{\begin{eqnarray}}
    \def\enl{\end{eqnarray}}
\def\bma{\begin{bmatrix}}
    \def\ema{\end{bmatrix}}
\def\bmx{\begin{matrix}}
    \def\emx{\end{matrix}}
\def\ben{\begin{enumerate}}
    \def\een{\end{enumerate}}
\def\bit{\begin{itemize}}
    \def\eit{\end{itemize}}
\def\bei{\begin{itemize}}
    \def\eei{\end{itemize}}
\def\bet{\begin{tabular}}
    \def\eet{\end{tabular}}

\newcommand{\ba}{\mathbf{a}}

\newcommand{\be}{\mathbf{e}}

\usepackage[dvipsnames]{xcolor}

\def\1{{\bf1}}

\def\bit{\begin{itemize}}
\def\eit{\end{itemize}}

\def\be{\begin{equation}}
\def\ee{\end{equation}}
\def\ba{\begin{eqnarray}}
\def\ea{\end{eqnarray}}

\def\bes{\begin{equation*}}
\def\ees{\end{equation*}}
\def\bas{\begin{eqnarray*}}
\def\eas{\end{eqnarray*}}
\def\as{{\it a.s.\ }}

\usepackage{url}
\usepackage{verbatim}

\newtheorem{Remark 1}{Remark}
\newtheorem{Remark 2}[Remark 1]{Remark}
\newtheorem{Remark 3}[Remark 1]{Remark}
\newtheorem{Remark 4}[Remark 1]{Remark}
\newtheorem{Remark 5}[Remark 1]{Remark}
\newtheorem{Remark 6}[Remark 1]{Remark}
\newtheorem{Remark 7}[Remark 1]{Remark}

\newtheorem{Lemma 1}{Lemma}
\newtheorem{Lemma 2}[Lemma 1]{Lemma}
\newtheorem{Lemma 3}[Lemma 1]{Lemma}
\newtheorem{Lemma 4}[Lemma 1]{Lemma}
\newtheorem{Lemma 5}[Lemma 1]{Lemma}
\newtheorem{Lemma 6}[Lemma 1]{Lemma}
\newtheorem{Lemma 7}[Lemma 1]{Lemma}

\newtheorem{Assumption 1}{Assumption}
\newtheorem{Assumption 2}[Assumption 1]{Assumption}
\newtheorem{Assumption 3}[Assumption 1]{Assumption}
\newtheorem{Assumption 4}[Assumption 1]{Assumption}
\newtheorem{Definition 1}{Definition}
\newtheorem{Theorem 1}{Theorem}
\newtheorem{Theorem 2}[Theorem 1]{Theorem}
\newtheorem{Theorem 3}[Theorem 1]{Theorem}
\newtheorem{Theorem 4}[Theorem 1]{Theorem}
\newtheorem{Theorem 5}[Theorem 1]{Theorem}
\newtheorem{Theorem 6}[Theorem 1]{Theorem}
\newtheorem{Theorem 7}[Theorem 1]{Theorem}
\newtheorem{Theorem 8}[Theorem 1]{Theorem}
\newtheorem{Theorem 9}[Theorem 1]{Theorem}
\newtheorem{Theorem 10}[Theorem 1]{Theorem}

\newtheorem{Proposition 1}{Proposition}

\title{\LARGE \bf
 Robust Constrained   Consensus and Inequality-constrained Distributed Optimization with Guaranteed Differential Privacy and Accurate Convergence}

\author{Yongqiang Wang, Angelia Nedi\'c
\thanks{The work was supported in part by the National Science Foundation under Grants   ECCS-1912702, CCF-2106293,  CCF-2106336, CCF-2215088, and CNS-2219487.}
\thanks{Yongqiang Wang is with the Department of Electrical and Computer Engineering, Clemson University, Clemson, SC 29634, USA
{\tt\small{yongqiw}@clemson.edu}
}%
\thanks{Angelia Nedi\'c is with the School of Electrical, Computer and Energy Engineering, Arizona State University, Tempe, AZ 85281, USA {\tt\small angelia.nedich@asu.edu}}
}

\begin{document}

\maketitle
\thispagestyle{empty}
\pagestyle{empty}

\begin{abstract}
We address differential privacy for fully distributed optimization subject to a shared inequality constraint. By co-designing the distributed optimization mechanism and the differential-privacy noise injection mechanism, we propose the first distributed constrained optimization algorithm that can ensure both provable convergence to a global optimal solution and rigorous $\epsilon$-differential privacy, even when the number of iterations tends to infinity. Our approach does not require the Lagrangian function to be strictly convex/concave, and allows the global objective function  to be non-separable. As a byproduct of the co-design, we also propose a new constrained consensus  algorithm that can achieve rigorous $\epsilon$-differential privacy while maintaining accurate convergence, which, to our knowledge,  has not been achieved before.   Numerical simulation results on a demand response control problem in smart grid  confirm the effectiveness of the proposed approach.
\end{abstract}

\section{Introduction}

In recent years,  distributed optimization subject to a shared inequality constraint is gaining increased traction due to its prevalence in power flow control \cite{chang2014distributed}, task  allocation and assignment \cite{notarnicola2019constraint}, cooperative model predictive control  \cite{patrascu2018convergence}, and  collective learning \cite{hershberger2001distributed}. Compared with the intensively studied   distributed optimization problem where   constraints do not  exist or  are  independent among agents, the incorporation of a shared inequality constraint   captures important practical  characteristics  such as  limited network resources or physical constraints \cite{notarstefano2019distributed}, and is mandatory in many applications of distributed optimization.

We consider the   distributed scenario where every agent  only has  access to its local cost  and constraint functions.  To account for the lack of information, agents exchange messages among local neighbors
to estimate  global information necessary for reaching a global optimal solution. To date, plenty of solutions have been proposed. The distributed dual subgradient method \cite{yang2010distributed,su2022convergence} is a commonly used approach to solving distributed optimization subject to a shared inequality constraint. However, this approach requires individual agents to globally solve local subproblems in every iteration, which is not  desirable when the local cost and constraint functions have some complex structure. To solve this issue, the distributed primal-dual subgradient method has been developed which only requires every   agent to simply evaluate the gradients of the cost and constraint functions in every iteration \cite{zhu2011distributed,yuan2011distributed}. This approach has also been extended in \cite{chang2014distributed} to relax the requirement for strict convexity/concavity of the Lagrangian function by introducing perturbations to the primal and dual updates.

  All of the aforementioned algorithms  require participating agents to repeatedly share iteration  variables, which are problematic when involved data are sensitive. For example, 
  in sensor network based target localization, the positions of sensors should  be kept private in sensitive (hostile) environments  ~\cite{zhang2019admm,huang2015differentially}. However, {in existing distributed optimization based localization   algorithms, the position  of a sensor is a parameter of its objective function, and as shown in ~\cite{zhang2019admm,huang2015differentially,burbano2019inferring}, it is easily inferable by an adversary  using  information shared in these distributed algorithms}. The privacy problem is more acute in distributed machine learning where involved training data may contain sensitive information such as medical or salary information (note that in machine learning, together with the model, training data  determines the objective function). In fact, as shown in our recent results \cite{wang2022quantization,wang2022tailoring,wang2022decentralized}, in the absence of  a   privacy mechanism, an adversary can use information shared  in distributed optimization to precisely recover the raw data used for training. It is worth noting that although distributed dual subgradient methods can avoid sharing the primal variable and only  share  dual variables  (see, e.g., \cite{falsone2017dual}), they  cannot provide strong privacy protection since the iteration trajectory of dual variables still bears information of the primal variable. This is particularly the case in primal-dual subgradient  methods like \cite{chang2014distributed,tjell2019privacy} and distributed dual  subgradient method \cite{han2021privacy} where disclosing dual updates in two consecutive iterations allows an observer to directly calculate some function value of the primal variable.

 Recently, plenty of efforts have been reported to address privacy protection in distributed optimization. One approach   is (partially) homomorphic encryption, which has been employed in both our own results ~\cite{zhang2019admm,zhang2018enabling}  and others ~\cite{freris2016distributed,lu2018privacy}. However, this approach leads to significant overheads in both communication and computation. Another approach is based on spatially or temporally correlated noises/uncertainties \cite{lou2017privacy,gade2018private,gao2022dynamics}, which can obfuscate information shared in distributed optimization. {This approach, unfortunately, has been shown in \cite{gao2022dynamics} to be susceptible to adversaries that are able to access all messages shared in the network}. Time-varying uncoordinated stepsizes \cite{wang2022decentralized,wang2022decentralized1} and stochastic quantization \cite{wang2022quantization} have also been exploited in our prior work to enable privacy protection in distributed optimization. With rigorous mathematical foundations yet   implementation simplicity and   post-processing immunity, Differential Privacy (DP)~\cite{dwork2014algorithmic} is   becoming a de facto standard for privacy protection.
  In fact, as DP  has  achieved remarkable successes in various
applications~\cite{han2016differentially,hale2017cloud,wang2017differential,zhang2019recycled,he2020differential},  efforts
have also emerged incorporating DP into distributed optimization.
  For example, DP-based privacy approaches have been proposed for distributed optimization by injecting DP noises
to exchanged messages~\cite{huang2015differentially,cortes2016differential,xiong2020privacy,ding2021differentially},
or objective functions~\cite{nozari2016differentially}. However, directly incorporating persistent DP noises into existing algorithms also unavoidably compromises the accuracy of optimization.   In addition, to our knowledge,  no results  have been reported able to achieve rigorous $\epsilon$-DP for both cost and constraint  functions in distributed  optimization subject to a shared inequality constraint. {Compared with privacy protection in constraint-free distributed optimization, privacy protection is more challenging in the presence of a shared inequality constraint: Not only does  the shared constraint  among agents introduce  an additional attack surface, convergence analysis also becomes  more involved under DP noises   because the entanglement of  amplitude-unbounded DP noises and nonlinearity induced by projection (necessary to address constraints) poses challenges to both optimality analysis and consensus characterization}.

In this paper, we propose a fully distributed constrained optimization algorithm that can {enable $\epsilon$-DP for the objective and constraint functions without losing accurate convergence}.
Our  basic idea is to suppress the influence of  DP noises by gradually weakening inter-agent interaction, which would  become unnecessary anyway
after convergence.   
Given that a key component of our  algorithm is constrained  consensus, we first propose a new robust constrained consensus  algorithm that {is robust to additive noises added to shared messages (for the purpose of DP) and can ensure accurate convergence even when these noises are persistent}.

The main contributions are summarized as follows:
\begin{enumerate}[wide, labelwidth=!,labelindent=8pt]
  \item  {We propose the first constrained consensus algorithm that can ensure both   accurate convergence and $\epsilon$-DP. This is in stark difference from others' results (e.g., \cite{huang2012differentially,nozari2015}) on differentially private consensus that have to sacrifice accurate convergence in order to ensure rigorous DP in the infinite horizon. Note that it is also different from our recent work in \cite{wang2022differentially_consensus} which does not consider constraints on the consensus variable. } 
  \item The proposed   fully distributed optimization algorithm, to  our knowledge,  is {the first  to achieve rigorous $\epsilon$-DP for both cost and constraint functions in distributed  constrained optimization}. Note that compared with constraint-free distributed optimization, the shared inequality constraint   introduces  an additional attack surface  and thus challenges in privacy protection. {In fact, to deal with the shared inequality constraint, a distributed optimization algorithm has to share  more variables among the agents, which leads to more noises injected into the algorithm (as every shared message has to be obfuscated by noise to enable DP). This makes the task of ensuring  a finite cumulative privacy budget in the infinite time horizon (under the constraint of accurate convergence) substantially more challenging than the constraint-free case.}
  \item  Besides achieving rigorous $\epsilon$-DP, the approach can simultaneously ensure  provable convergence to an optimal solution,  which is in sharp contrast to existing DP solutions for constraint-free distributed optimization and games (e.g., \cite{huang2015differentially,ye2021differentially,wang2022differentially}) that have to trade convergence accuracy for DP.
   \item By employing the primal-dual perturbation technique, our  algorithm can ensure convergence to a global optimal solution without requiring the Lagrangian function to be strictly convex/concave.
     {Note that different from~\cite{chang2014distributed} which constructs optimal solutions using  running averages of its iterates, we  prove that our iterates    directly converge to an optimal solution, which enables us to avoid the additional computational and storage overheads  for maintaining additional running averages of iterates in \cite{chang2014distributed}. In addition, to deal with DP noises, our algorithm structure and parameters are significantly different from \cite{chang2014distributed} {which} does not consider DP noises.}
       Our formulation also allows the global objective function  to be non-separable, which is more general than the intensively studied  distributed optimization problem with separable objective functions.
  \item  Even without considering privacy, our proof techniques are  of independent interest. In fact, to our knowledge, our convergence analysis is the first  to prove accurate convergence under entangled dynamics of unbounded DP noises and projection-induced nonlinearity in distributed optimization.
  \end{enumerate}

The organization of the paper is as follows.
Sec.~\ref{sec-problem} gives the problem formulation
and some preliminary results. Sec.~\ref{se:algorithm1} presents a new constrained consensus  algorithm that is necessary for the development of our distributed  optimization algorithm. This section also proves that  the algorithm can ensure both almost sure convergence and rigorous $\epsilon$-DP with a finite cumulative privacy budget, even when the number of iterations tends to infinity. Sec.~\ref{se:algorithm_2} proposes a fully distributed constrained optimization algorithm, while its  almost sure convergence to a global optimal solution   and rigorous $\epsilon$-DP guarantee are presented in Sec. \ref{se:convergence_algo2} and Sec.~\ref{se:DP_algo2}, respectively.  Sec.~\ref{se:simulation} presents numerical simulation results to confirm the obtained results. Finally, Sec.~\ref{se:conclusions} concludes the paper.

{\bf Notations:}
 We use $\mathbb{R}^d$ to denote the Euclidean space of
dimension $d$.  
We write $I_d$ for the identity matrix of dimension $d$,
and ${\bf 1}_d$ for  the $d$-dimensional  column vector with all
entries equal to 1; in both cases  we suppress the dimension when it is
clear from the context. 
  We use $\langle\cdot,\cdot\rangle$ to denote the inner product and
 $\|x\|$ for the standard Euclidean norm of a vector $x$. We use $\|x\|_1$ for the
 $\ell_1$ norm of a vector $x$.
We write $\|A\|$ for the matrix norm induced by the vector norm $\|\cdot\|$.
 $A^T$ denotes the transpose of a matrix $A$.
 Given vectors $x_1,\cdots,x_m$, we define ${\rm col}(x_1,\cdots,x_m)=[x_1^T,\cdots,x_m^T]^T$. We use $\otimes$ to denote the Kronecker product.
  Often, we abbreviate {\it almost surely} by {\it a.s}. {We use the superscript $k$ on the shoulder of a symbol to denote its value in the $k$th iteration  throughout the paper.}


\section{Problem Formulation and Preliminaries}\label{sec-problem}
 \subsection{Problem Formulation}
We consider a distributed optimization problem  among a set of $m$  agents  $[m]=\{1,\,\ldots,m\}$.
Agent $i$ is characterized by a local optimization variable $x_i\in \mathbb{R}^d$, a local constraint set   $\mathcal{X}_i\subseteq \mathbb{R}^d$, and a local mapping
$f_i:\mathbb{R}^d\rightarrow \mathbb{R}^M$.
The agents work together to minimize a network objective function
\[
J(x_1,\ldots,x_m)\triangleq F\big({  \textstyle\sum_{i=1}^{m}}f_i(x_i)\big),
\]
where $F:\mathbb{R}^M\rightarrow \mathbb{R} $ and $J:\mathbb{R}^{md}\rightarrow \mathbb{R} $.

Moreover, the optimization variables of all agents  must satisfy a shared  inequality constraint
\begin{equation}{\label{eq:constraint}}
\textstyle\sum_{i=1}^{m} g_i(x_i)\leq 0,
\end{equation}
 where $g_i:\mathbb{R}^d\rightarrow \mathbb{R}^P$ is a mapping with components $g_{ij}: \mathbb{R}^d\rightarrow \mathbb{R}$ for $j=1,\ldots,P$ , i.e., $g_i({x_i})=[g_{i1}({x_i}),\ldots,g_{iP}({x_i})]^T$.   The vector inequality in~(\ref{eq:constraint}) means that each element of the vector $\sum_{i=1}^{m} g_i(x_i)$ is nonpositive.
 Such  inequality constraints arise in various  scenarios involving upper or lower limits of shared resources, and have been a common assumption in smart grid control \cite{chang2014distributed,camisa2019distributed}, distributed model predictive control \cite{necoara2013rate}, and distributed sparse regression problems \cite{chang2013distributed}, etc. It is worth noting that linear equality  constraints can also be dealt with under our framework by means of double-sided inequalities.
Compared with the commonly studied scenario where the objective functions are  additively separable across agents \cite{nedic2009distributed}, here we allow the global objective function to be non-separable (described by $F(\cdot)$), which is more general and necessary in  some  applications such as  flow control and  regression problems \cite{chang2014distributed,sundhar2012new}.

The distributed optimization problem  can be formalized as follows:
\begin{equation}\label{eq:formulation}
\begin{aligned}
&\qquad \min_{x_1,\ldots,x_m}J(x_1,\ldots,x_m)\\
&  {\rm s.t.}\quad x_i\in\mathcal{X}_i\:\: {\rm and}\:\:  \textstyle\sum_{i=1}^{m}g_i(x_i)\leq 0.
\end{aligned}
\end{equation}

 We use the following assumption   analogous to that of~\cite{chang2014distributed}, with minor adjustments to take into account that $f_i$ and $g_i$ are mappings.
\begin{Assumption 1}\label{ass:compact}
For all $i\in[m]$, the following conditions hold:
 \begin{itemize}
 \item[(a)]
 The local constraint set $\mathcal{X}_i$ is nonempty, compact, and convex.
 \item[(b)]
 The mapping $f_i(\cdot)$
 has continuous Jacobian $\nabla f_i(\cdot)$ over an open set containing the constraint set $\mathcal{X}_i$.
 \item[(c)]
 The component function $g_{ij}(\cdot):\mathbb{R}^d\to\mathbb{R}$ is convex for all $j=1,\ldots,P$, and
 the constraint mapping $g_i(\cdot)$ has Lipschitz continuous Jacobian $\nabla g_i(\cdot)$ over an open set containing the constraint set $\mathcal{X}_i$.
 \end{itemize}
\end{Assumption 1}

\def\cX{{\mathcal{X}}}
We let $D_{\cX}>0$ be a constant such that
 $\|x_i\|\leq D_\cX  $ for all $x_i\in\cX_i$.
We next discuss the implications of Assumption~\ref{ass:compact}.
Since the set $\cX_i$ is compact and the Jacobian $\nabla f_i(\cdot)$ is continuous over $\cX_i$,
it follows that $\nabla f_i(\cdot)$ is bounded over $\cX_i$. Consequently, the mapping $f_i(\cdot)$
is Lipschitz continuous on $\cX_i$. Therefore, there exists a constant $L_f$ such that the following relations hold for all $i\in[m]$
and all $u,v\in \mathcal{X}_i$:
  \[
  \begin{aligned}\|\nabla f_i(u)\|&\leq L_f, \quad
\|f_i(u)-f_i(v)\| \leq L_f\left\| u-v\right\|.
  \end{aligned}
  \]

 Similarly, since the Jacobian $\nabla g_i(\cdot)$ is continuous over the compact set $\cX_i$,
the Jacobian $\nabla g_i(\cdot)$ is bounded on $\cX_i$ and
the mapping $g_i(\cdot)$ is Lipschitz continuous on $\cX_i$. These statements and the assumption that $\nabla g_i(\cdot)$ is
Lipschitz continuous on $\cX_i$ imply that there exist constants $L_g,G_g>0$ such that
the following relations are valid for all $u,v\in \mathcal{X}_i$:
  \[
  \begin{aligned}
   &\|\nabla g_i(u) \|\leq L_g, \qquad
  \|g_i(u)-g_i(v)\| \leq L_g\left\| u-v\right\|,\\
  &\|\nabla g_i(u)-\nabla g_i(v) \|\leq G_g\|u-v\|.
  \end{aligned}
  \]

 Consequently, the (Lipschitz) continuity of $f_i(\cdot)$ and $g_i(\cdot)$
imply that both $f_i(\cdot)$ and $g_i(\cdot)$
are bounded on the compact set $\cX_i$.
Hence, there exist constants $C_f,C_g>0$ such that for all $i\in[m]$ and all $x_i\in\mathcal{X}_i$:
\begin{equation}\label{eq:C_fC_g}
\|f_i(x_i)\| \leq C_f,\quad \|g_i(x_i)\|\leq C_g.
\end{equation}

We   make the following assumption on  $F(\cdot)$ and $J(\cdot)$:
\begin{Assumption 2}\label{ass:network_function_Lipschitz}
 The function $F(\cdot)$ satisfies the following  relationship for all  $v,u\in\mathbb{R}^M$:
  \[
  \begin{aligned}
\|\nabla F(v)-\nabla F(u)\| \leq G_F\left\| v-u\right\|,  \quad   \|\nabla F(v)\| \leq L_F.
  \end{aligned}
  \]
  The function $J(\cdot)$ is convex and satisfies the following  relationship for all  $x,y\in\mathcal{X}\triangleq \mathcal{X}_1\times\cdots\times\mathcal{X}_m$:
  \[
  \begin{aligned}
\|\nabla J(x)\hspace{-0.07cm}-\hspace{-0.07cm}\nabla J(y)\|\leq G_J\hspace{-0.07cm}\left\| x\hspace{-0.07cm}-\hspace{-0.07cm}y\right\|,\:   \| J(x)\hspace{-0.07cm}-\hspace{-0.07cm}J(y)\| \leq L_J\hspace{-0.03cm}\|\hspace{-0.03cm}x-\hspace{-0.07cm}y\|.
  \end{aligned}
  \]
\end{Assumption 2}


We consider the scenario in which  each agent $i$ can access $F(\cdot)$, $f_i(\cdot)$, $g_i(\cdot)$, and $\mathcal{X}_i$ only. Hence, agents have to share information with their immediate neighbors to cooperatively minimize the network cost $J(\cdot)$.
 We describe the local interaction using a weight matrix
$W=\{w_{ij}\}$, where $w_{ij}>0$ if agent  $j$ and agent $i$ can directly communicate with each other,
and $w_{ij}=0$ otherwise. For an agent $i\in[m]$,
its  neighbor set
$\mathbb{N}_i$ is defined as the collection of agents $j$ such that $w_{ij}>0$.
We define $w_{ii}\triangleq-\sum_{j\in\mathbb{N}_i}w_{ij}$  for all $i\in [m]$. Furthermore,
We make the following assumption on the matrix $W$.

\begin{Assumption 1}\label{as:L}
 The matrix  $W=\{w_{ij}\}\in \mathbb{R}^{m\times m}$ is symmetric and satisfies
    ${\bf 1}^TW={
  0}^T$, $W{\bf 1}={
  0}$, and $ \|I+W-\frac{{\bf 1}{\bf 1}^T}{m}\|<1$.
\end{Assumption 1}

Assumption~\ref{as:L} ensures  that the graph induced by the matrix $W$ is connected, i.e.,
there is a  path in the graph
from each agent to every other agent.

\subsection{Primal-Dual Method}

In primal-dual based approach we consider the Lagrangian dual problem of (\ref{eq:formulation}):
\begin{equation}\label{eq:Lagrange_dual}
  \max_{\lambda}\{\min_{x\in\mathcal{X}}\mathcal{L}(x,\lambda) \},
\end{equation}
where $\lambda\in\mathbb{R}^P_+$ (the non-negative orthant of $\mathbb{R}^P$) is the dual variable associated with the constraint $\sum_{i=1}^{m}g_i(x_i)\leq 0$, and $\mathcal{L}(x,\lambda)$ is the Lagrangian function defined as $\mathcal{L}(x,\lambda)=J(x_1,\ldots,x_m)+\lambda^T\sum_{i=1}^{m}g_i(x_i)$.
Under Assumption \ref{ass:compact} and Assumption \ref{ass:network_function_Lipschitz}, $\mathcal{L}(x,\lambda)$ is convex in $x$. We assume that the Slater condition holds, and hence, the strong duality holds for (\ref{eq:formulation}), and problem (\ref{eq:formulation}) can be solved by considering its dual problem (\ref{eq:Lagrange_dual}) (see, for example~\cite{boyd2004convex,bno2003}).
A classical approach to problem (\ref{eq:Lagrange_dual}) is the dual subgradient method \cite{boyd2004convex,bno2003}.
This method  requires having a solution for the inner minimization problem in (\ref{eq:Lagrange_dual}) at each iteration, which makes the iterate updates computationally expensive.  For this reason,
we employ the primal-dual gradient method that takes a single step for the minimization problem \cite{uzawa1958iterative}. More specifically, in each iteration $k$, two updates are performed as follows:
\begin{equation}\label{eq:PD_centralized}
  \begin{aligned}
         x^{k+1}&=\Pi_{\mathcal{X}}\left[x^k-\gamma^k  \mathcal{L}_x(x^k,\lambda^k)\right],\\
         \lambda^{k+1}&=\Pi_{\mathbb{R}^P_+}\left[\lambda^k+\gamma^k \mathcal{L}_\lambda(x^k,\lambda^k) \right],
  \end{aligned}
\end{equation}
where $\Pi_{\cX}[\cdot]$ denotes Euclidean projection of a vector   on the set $\cX$, $\gamma^k > 0$ is the stepsize, and  $\mathcal{L}_x(x^k,\lambda^k)$ and $\mathcal{L}_\lambda(x^k,\lambda^k)$ denote the gradient of $\mathcal{L}(x^k,\lambda^k)$ with respect to $x$ and $\lambda$, respectively, i.e.,
\begin{equation}\label{eq:L_x}
\begin{aligned}
\mathcal{L}_x(x^k,\lambda^k)\triangleq {\rm col}(\mathcal{L}_{x_1}(x^k,\lambda^k),\ldots,\mathcal{L}_{x_m}(x^k,\lambda^k)),
\end{aligned}
\end{equation}
with
\begin{equation}\label{eq:L_x_i}
 \hspace{-0.1cm}\mathcal{L}_{x_i}(x^k,\lambda^k)\triangleq \nabla\hspace{-0.05cm} f_i^T(x_i^k)\hspace{-0.03cm}\nabla F\hspace{-0.1cm}\left({\textstyle\sum_{j=1}^{m}}f_j(x_j^k)\right)\hspace{-0.04cm}+\hspace{-0.04cm}\nabla g_i^T(x_i^k)\lambda^k
\end{equation}
for $i\in[m]$, and
\begin{equation}\label{eq:L_lambda}
   \mathcal{L}_{\lambda}(x^k,\lambda^k)\triangleq {\textstyle\sum_{i=1}^{m}}g_i(x_i^k).
\end{equation}

It is well known that if the primal-dual iterates of~(\ref{eq:PD_centralized}) converge to a saddle point  $ (x^\ast,\lambda^\ast)$, i.e., $\mathcal{L}(x^\ast,\lambda)\leq \mathcal{L}(x^\ast,\lambda^\ast)\leq \mathcal{L}(x,\lambda^\ast)$ holds for any $x\in\mathcal{X}$ and $\lambda\geq 0$, then the obtained $x^\ast$ is an optimal solution to the original problem (\ref{eq:formulation})
(see~\cite{bno2003}).
However, to ensure convergence of the sequence $(x^k,\lambda^k)$ to a saddle point, it is often assumed that the Lagrangian function $\mathcal{L}(x,\lambda)$ is strictly convex in $x$ and strictly concave in $\lambda$, which, however, does not hold  here since the Lagrangian function $\mathcal{L}(x,\lambda)$ is linear in $\lambda$.
Therefore, to circumvent this problem, inspired by \cite{chang2014distributed}, we use the primal-dual perturbation gradient method. The basic idea  is to first obtain some  perturbed versions of $x^k$ and $\lambda^k$ (represented as $ \alpha^k$ and $ {\beta}^k$):
 \begin{equation}\label{eq:PD_centralized_per}
  \begin{aligned}
         \alpha^{k+1}&=\Pi_{\mathcal{X}}\left[x^k-\rho_1  \mathcal{L}_x(x^k,\lambda^k)\right],\\
         \beta^{k+1}&=\Pi_{\mathbb{R}^P_+}\left[\lambda^k+\rho_2 \mathcal{L}_\lambda(x^k,\lambda^k) \right],
  \end{aligned}
\end{equation}
where $\rho_1>0$ and $\rho_2>0$ are some constants.
Then, the approach updates  the primal and dual variables as follows:
\begin{equation}\label{eq:PD_centralized_perburb}
  \begin{aligned}
         x^{k+1}&=\Pi_{\mathcal{X}}\left[x^k-\gamma^k  \mathcal{L}_x(x^k, \beta^{k+1})\right],\\
         \lambda^{k+1}&=\Pi_{\mathbb{R}^P_+}\left[\lambda^k+\gamma^k \mathcal{L}_\lambda( \alpha^{k+1},\lambda^k) \right].
  \end{aligned}
\end{equation}
It has been proven in \cite{kallio1994perturbation,kallio1999large} that by carefully designing the stepsize, this primal-dual perturbation gradient method can ensure convergence to a saddle point of $\mathcal{L}(x,\lambda)$ without imposing any strict convexity/concavity conditions on  $\mathcal{L}(x,\lambda)$.

\subsection{Preliminaries}

We   need  the following lemmas for convergence analysis:

\begin{Lemma 1}[\cite{polyak87}, Lemma 11, page 50]\label{Lemma-polyak1}
Let $\{v^k\}$, $\{u^k\}$, $\{a^k \}$, and $\{b^k\}$ be random nonnegative scalar sequences such that
$\sum_{k=0}^\infty { a^k}<\infty$ and $\sum_{k=0}^\infty {b^k}<\infty$ \as
and
\[
\begin{aligned}
\mathbb{E}\left[v^{k+1}|\mathcal{F}^k \right]\le(1+{a^k}) v^k -u^k+{ b^k},\quad \forall k\geq 0 \quad \as
\end{aligned}
\]
where 
$\mathcal{F}^k=\{v^\ell,u^\ell,{a^\ell},{b^\ell};\, 0\le \ell\le k\}$.
Then $\sum_{k=0}^{\infty}u^k<\infty$ and $\lim_{k\to\infty}v^k=v$ hold {\it a.s.} for a random variable $v\geq 0$.
\end{Lemma 1}

\begin{Lemma 2}\cite{wang2022tailoring}\label{Lemma-polyak}
Let $\{v^k\}$,$\{{a^k}\}$, and $\{p^k\}$ be random nonnegative scalar sequences, and
$\{q^k\}$ be a deterministic nonnegative scalar sequence satisfying
$\sum_{k=0}^\infty {a^k}<\infty$  {\it a.s.},
$\sum_{k=0}^\infty q^k=\infty$, $\sum_{k=0}^\infty p^k<\infty$ {\it a.s.},
and
\[
\mathbb{E}\left[v^{k+1}|\mathcal{F}^k\right]\le(1+{ a^k}-q^k) v^k +p^k,\quad \forall k\geq 0\quad\as
\]
where $\mathcal{F}^k=\{v^\ell,{a^\ell},p^\ell; 0\le \ell\le k\}$.
Then, $\sum_{k=0}^{\infty}q^k v^k<\infty$ and
$\lim_{k\to\infty} v^k=0$ hold almost surely.
\end{Lemma 2}

\subsection{On Differential Privacy}\label{se:definition}
Since  iterative optimization algorithms  involve continual information sharing among   agents, we use the notion of $\epsilon$-DP for continuous bit streams \cite{dwork2010differential}, which has recently been applied to constraint-free distributed optimization  (see \cite{huang2015differentially} as well as our own work \cite{wang2022quantization,wang2022tailoring}). {We consider DP for the objective and constraint functions $f_i$ and $g_i$ of individual agents.} To enable $\epsilon$-DP, we  add Laplace noise to all shared messages {(note that similar results can be obtained if Gaussian noise is used)}. For a constant $\nu>0$, we use ${\rm Lap}(\nu)$ to denote a Laplace distribution    with the probability density function $x\mapsto\frac{1}{2\nu}e^{-\frac{|x|}{\nu}}$.  ${\rm Lap}(\nu)$'s mean is   zero   and its  variance is $2\nu^2$.
To facilitate DP analysis, we represent the distributed  optimization problem in (\ref{eq:formulation}) by two {letters} ($\mathcal{J},\mathcal{G}$), where
    $ \mathcal{J} \triangleq\{ f_1,\,\ldots,f_m\}$ is the set of mapping functions of individual agents, and $ \mathcal{G} \triangleq\{ g_1,\,\ldots,g_m\}$ is the set of constraint functions of individual agents. Then, we define adjacency between distributed optimization problems  as follows:

\begin{Definition 1}\label{de:adjacency}
Two distributed constrained optimization problems {$\mathcal{P}=(\mathcal{J}, \mathcal{G})$} and {$ \mathcal{P}'=(\mathcal{J}', \mathcal{G}')$} are adjacent if the following conditions hold:
\begin{itemize}
\item there exists an $i\in[m]$ such that ${\rm col}(f_i,g_i)\neq {\rm col}(f_i',g_i')$ but ${\rm col}(f_j,g_j)= {\rm col}(f_j',g_j')$ for all $j\in[m],\,j\neq i$;
\item   ${\rm col}(f_i,g_i)$ and ${\rm col}(f_i',g_i')$, which are not the same,  have similar behaviors   around $\theta^\ast$, the  solution of $\mathcal{P}$. More specifically,  there exits some $\delta>0$  such that for all $v$ and $v'$ in  $B_\delta(\theta^\ast)\triangleq\{u:u\in\mathbb{R}^d, \|u-\theta^\ast\|<\delta\}$, we have $f_i(v)= f'_i(v')$ and $g_i(v)=g'_i(v')$.
\end{itemize}
\end{Definition 1}

 \begin{Remark 1}
 In Definition \ref{de:adjacency}, since  the difference between ${\rm col}(f_i,g_i)$ and ${\rm col}(f_i',g_i')$ in the first condition  can be arbitrary,
 additional restrictions have to be imposed to ensure rigorous DP in distributed optimization.  {Different from \cite{huang2015differentially,chen2023differentially} which restrict  all gradients  to be uniformly bounded, and \cite{huang2024differential} which restricts the gradient difference to be uniformly bounded,}  we add the {second} condition, which, as shown later, {together with the proposed noise-robust algorithm,} allows us to ensure  rigorous DP while maintaining provable convergence to the optimal solution. {The condition requires that  two adjacent problems have the same solutions.} It  can be satisfied by a large class of functions that do not necessarily satisfy the  bounded-gradient condition or the bounded-gradient-difference condition.
   For example,  any two functions in the form of $h_1(v)\cdot h_2(v)$ can be $f_i$ and $f'_i$, where $h_1(v)$ is any {nonnegative} smooth  function and $h_2(v)$ is given by
 \[
 h_2(v)=\left\{\begin{array}{cl} 0& \|v-\theta^\ast\| < \delta\\
                               (\|v-\theta^\ast\|-\delta)^2 &\|v-\theta^\ast\| \geq \delta

                                      \end{array}
                                       \right.
  \]
 Note that $g_i$ and $g'_i$ can be constructed in the same manner.
\end{Remark 1}
\begin{Remark 1}
It is worth noting that compared with  constraint-free distributed optimization, in distributed optimization subject to a shared inequality constraint, we have to take  care of the privacy of constraint functions besides individual cost functions. In fact, since the constraint functions are coupled among the agents, they introduce an additional  attack surface and make privacy protection more challenging.
\end{Remark 1}

 Given a distributed constrained optimization algorithm, we represent an execution of such an algorithm as $\mathcal{A}$, which is an infinite sequence of the   iteration variable $\vartheta$, i.e., $\mathcal{A}=\{\vartheta^0,\vartheta^1,\ldots\}$ (note that $\vartheta$ can be an augmentation  of several variables). We consider adversaries able to observe all shared messages among the agents. Therefore, the observation part of an execution is the infinite sequence of shared messages, which is denoted as $\mathcal{O}$. Given a distributed optimization problem $\mathcal{P}$ and an initial state $\vartheta^0$, we define the observation mapping as $\mathcal{R}_{\mathcal{P},\vartheta^0}(\mathcal{A})\triangleq \mathcal{O}$. Given a distributed   optimization  problem $\mathcal{P}$, observation sequence $\mathcal{O}$, and an initial state $\vartheta^0$,  $\mathcal{R}_{\mathcal{P},\vartheta^0}^{-1}(\mathcal{O})$ is the set of executions $\mathcal{A}$ that can generate the observation $\mathcal{O}$.
 \begin{Definition 1}
   ($\epsilon$-differential privacy, adapted from \cite{huang2015differentially}). For a given $\epsilon>0$, an iterative  algorithm {$\mathcal{A}$} is $\epsilon$-differentially private if for any two adjacent $\mathcal{P}$ and $\mathcal{P}'$, any set of observation sequences $\mathcal{O}_s\subseteq\mathbb{O}$ (with $\mathbb{O}$ denoting the set of all possible observation sequences), and any initial state ${\vartheta}^0$, the following relationship always holds
    \begin{equation}
        \mathbb{P}[\mathcal{R}_{\mathcal{P},\vartheta^0}{(\mathcal{A})}\in \mathcal{O}_s]\leq e^\epsilon\mathbb{P}[\mathcal{R}_{\mathcal{P}',\vartheta^0}{(\mathcal{A})}\in \mathcal{O}_s],
    \end{equation}
    with the probability $\mathbb{P}$  taken over the randomness over iteration processes.
 \end{Definition 1}

The adopted $\epsilon$-DP definition ensures that an adversary having access to all shared messages cannot gain information with a  significant probability of any participating agent's mapping or constraint function. One can see that a smaller $\epsilon$ means a higher level of privacy protection. 


\section{Differentially-private constrained consensus}\label{se:algorithm1}
A necessary building block of our differentially-private distributed constrained optimization algorithm is the constrained consensus  technique. However,  the conventional constrained consensus  algorithm (see, e.g., \cite{nedic2010constrained}) cannot be used since   added DP noises may  stop the agents from reaching consensus.
Hence, in this section,   we first propose a differentially-private constrained consensus algorithm that can achieve both provable consensus and rigorous $\epsilon$-DP for the case where the agents' decisions are restricted to a {\it convex closed set $X\subseteq\mathbb{R}^d$}.
   To our knowledge, this is the first time that both goals are achieved in constrained consensus. The key idea of our algorithm is to use a weakening factor (represented as $\chi^k$) to gradually eliminate the influence of DP noises on the consensus accuracy. It is worth noting that different from the conventional constraint-free average consensus problem which has a unique final consensus value,    constrained consensus has many possible converging points. Therefore, in practical applications, {usually an additional input signal  is injected by each agent to ensure convergence to a desired value. For example, in distributed optimization, the gradient of each agent is used as this input signal to ensure convergence to an optimal solution \cite{nedic2010constrained} (also see the update of dual variables in (\ref{eq:primal_dual_update}) of Algorithm 2 in Sec. IV.)}.  Therefore, we let each agent apply an input signal (represented as $r_i$) to control the final consensus value. The  algorithm is summarized below.

\noindent\rule{0.49\textwidth}{0.5pt}
\noindent\textbf{Algorithm 1: Differentially-private constrained consensus}
\noindent\rule{0.49\textwidth}{0.5pt}
\begin{enumerate}[wide, labelwidth=!, labelindent=0pt]
    \item[] Parameters: {Deterministic } sequences $\chi^k>0$ and  $\gamma^k>0$.
    \item[] Every agent $i$'s  input is $r_i^k$.
    Every agent $i$ maintains one state variable  $x_i^k$, which is initialized randomly in $X$.
    \item[] {\bf for  $k=1,2,\ldots$ do}
    \begin{enumerate}
        \item Every agent $j$ adds  DP noise   $\zeta_j^{k}$ 
        to its state
    $x_j^k$,  and then sends the obscured state ${y_j^k\triangleq}x_j^k+\zeta_j^{k}$ to agent
        $i\in\mathbb{N}_j$.
        \item After receiving  ${y_j^k}$ from all $j\in\mathbb{N}_i$, agent $i$ updates its state  as follows:
        \begin{equation}\label{eq:update_in_Algorithm1}
\begin{aligned}
            \hspace{-0.3cm} x_i^{k+1}
              =\Pi_X\left[ x_i^k+\chi^k{ \textstyle\sum_{j\in \mathbb{N}_i}} w_{ij}({y_j^k}-x_i^k)+\gamma^k r_i^k\right].
        \end{aligned}
        \end{equation}
    \end{enumerate}
\end{enumerate}
\vspace{-0.1cm} \rule{0.49\textwidth}{0.5pt}
The set $X\subset\mathbb{R}^d$ is assumed to be nonempty, closed, and convex.
We make the following assumption on the DP noises:
\begin{Assumption 1}\label{ass:dp-noise}
For every $i\in[m]$ and every $k$, conditional on the state $x_i^k$,
the  DP noise $\zeta_i^k$ satisfies $
\mathbb{E}\left[\zeta_i^k\mid x_i^k\right]=0$ and
$\mathbb{E}\left[\|\zeta_i^k\|^2\mid x_i^k\right]=(\sigma_{i}^k)^2$ for all  $k\ge0$, and
\begin{equation}\label{eq:condition_assumption1}
\sum_{k=0}^\infty (\chi^k)^2\, \max_{i\in[m]}(\sigma_{i}^k)^2 <\infty,
\end{equation} where $\{\chi^k\}$ is the weakening factor sequence from Algorithm~1. Moreover,
$\mathbb{E}\left[\|x_i^0\|^2\right]<\infty$ holds for  $\forall i\in[m]$.
\end{Assumption 1}
\begin{Remark 1}
Note that Assumption \ref{ass:dp-noise} allows the variance of noise to increase with time. For example, when $\chi^k$ is set as $\mathcal{O}(\frac{1}{k^{0.7}})$, the condition in (\ref{eq:condition_assumption1}) can be guaranteed even if $\sigma_i^k$ increases with time, as long as its increasing rate is no faster than $\mathcal{O}(k^{0.2})$. In fact, as elaborated later on, allowing $\sigma_i^k$ to increase with time is key for our approach to enable both accurate convergence and rigorous $\epsilon$-DP.
\end{Remark 1}

\subsection{Convergence Analysis}

\begin{Lemma 1}\label{Lemma:projection} (Lemma 1 of \cite{nedic2010constrained})
 Let $X$ be a nonempty closed convex set in $\mathbb{R}^d$, then for any $x\in\mathbb{R}^d$ and $y\in X$, we have
 \[
 \|\Pi_{X}[x]-y\|\leq \|x-y\|^2-\| \Pi_{X}[x]-x \|^2.
 \]
\end{Lemma 1}

\begin{Lemma 1}\label{Le:rho_2} (Lemma 4 of \cite{wang2022differentially_consensus})
Under Assumption \ref{as:L}, for a positive sequence  $\{\chi^k\}$  satisfying $
\sum_{k=0}^\infty (\chi^k)^2<\infty
$, there  exists  a $T\geq 0$ such that
 $
\|I- \frac{{\bf 1}{\bf 1}^T}{m}+\chi^k W\|< 1-\chi^k |{\rho_2}|
$
holds for all $k\geq T$, where ${\rho_2}$ is $W$'s second largest eigenvalue.
\end{Lemma 1}


{
\begin{Assumption 1}\label{ass:bounded_r}
 There exists a constant $C$  such that $\|r_i^k\|<C$ holds for all $i\in [m]$ and $k\geq 0$.
\end{Assumption 1}

Defining $\bar{x}^k=\frac{\sum_{i=1}^{m}x_i^k}{m}$, we have the following results:
}

\begin{Theorem 1}\label{Th:consensus_tracking}
 Under {Assumptions  \ref{as:L}, \ref{ass:dp-noise}, and \ref{ass:bounded_r},} if
   the nonnegative sequences $\{\gamma^k\}$  and $\{\chi^k\}$ in Algorithm~1   satisfy
    \begin{equation}\label{eq:condtions_chi}
     \sum_{k=0}^{\infty}\chi^k=\infty,\, \sum_{k=0}^{\infty}(\chi^k)^2<\infty,\, \sum_{k=0}^{\infty}\frac{(\gamma^k)^2}{\chi^k}<\infty,
    \end{equation}
     then, the following results hold almost surely:
     \begin{enumerate}
    \item   $\lim_{k\to\infty}\|x_i^k - \bar{x}^k\|=0$  for all $i\in[m]$.
    \item    $\sum_{k=0}^{\infty}\chi^k\sum_{i=1}^{m}\|x_i^k-\bar{x}^k\|^2<\infty$;
    \item   $\sum_{k=0}^{\infty}\gamma^k \sum_{i=1}^{m}\|x_i^k-\bar{x}^k\| <\infty$.
  \end{enumerate}
\end{Theorem 1}

\begin{proof}
According to Lemma \ref{Lemma:projection} {on page 5} and the update rule of $x_i^k$ in (\ref{eq:update_in_Algorithm1}), we  have the following relationship for all $i\in[m]$ and   $y\in X$:
\[
 \|x_i^{k+1}-y\|^2\hspace{-0.06cm}\leq\hspace{-0.06cm} \left\|x_i^k\hspace{-0.06cm}+\hspace{-0.06cm}\chi^k{\textstyle\sum_{j\in \mathbb{N}_i}} w_{ij}(x_j^k+\zeta_j^k\hspace{-0.06cm}-\hspace{-0.06cm}x_i^k)\hspace{-0.06cm}+\hspace{-0.06cm}\gamma^k r_i^k\hspace{-0.06cm}-\hspace{-0.06cm}y\right\|^2 .
\]

Plugging $y=\bar{x}^k\in X$ to the preceding inequality, we  obtain
 \[
 \|x_i^{k+1}\hspace{-0.07cm}-\bar{x}^k\|^2\hspace{-0.08cm}\leq\hspace{-0.08cm} \left\|x_i^k\hspace{-0.06cm}+\hspace{-0.06cm}\chi^k{\textstyle\sum_{j\in \mathbb{N}_i}}\hspace{-0.06cm} w_{ij}(x_j^k\hspace{-0.06cm}+\hspace{-0.06cm}\zeta_j^k\hspace{-0.06cm}-\hspace{-0.06cm}x_i^k)\hspace{-0.06cm}+\hspace{-0.06cm}\gamma^k r_i^k\hspace{-0.07cm}-\hspace{-0.06cm}\bar{x}^k\right\|^2 .
\]
Given $\sum_{i=1}^{m}\|x_i^{k+1}-\bar{x}^k\|^2=\|x^{k+1}-{\bf 1}\otimes \bar{x}^k\|^2$, with $\otimes$ denoting the Kronecker product, we can add both sides of the preceding inequality from $i=1$ to $i=m$ to obtain
\[
\|x^{k+1}\hspace{-0.05cm}-\hspace{-0.05cm}{\bf 1}\otimes \bar{x}^k\|^2\hspace{-0.05cm}\leq\hspace{-0.05cm} \| ((I+\chi^k W)\otimes I_d ) x^k+\chi^k\zeta_w^k+\gamma^kr^k\hspace{-0.05cm}-\hspace{-0.05cm}{\bf 1} \otimes \bar{x}^k \|^2,
\]
where  $x^k={\rm col}(x_1^k,\ldots, x_m^k)$, $r^k={\rm col}(r_1^k,\ldots, r_m^k)$ and
$\zeta_w^k={\rm col}(\zeta_{w1}^k,\ldots,\zeta_{wm}^k)$ with $\zeta_{wi}^k\triangleq \sum_{j=1,j\neq i}^{m}w_{ij}\zeta_j^k$.

 By defining $W^k\triangleq  (I+\chi^k W-\frac{{\bf 1} {\bf 1}^T}{m})\otimes I_d$ and using
$\bar{x}^k=\frac{{\bf 1}^T\otimes I_d}{m}x^k$, we can further simplify the preceding inequality as
\begin{equation}\label{eq:consensus1}
\begin{aligned}
&\|x^{k+1}-{\bf 1}\otimes \bar{x}^k\|^2\leq \| W^k x^k+\chi^k\zeta_w^k+\gamma^kr^k \|^2\\
&=\| W^k x^k+\gamma^kr^k \|^2+\|\chi^k\zeta_w^k\|^2+2\langle W^k x^k+\gamma^kr^k,  \chi^k\zeta_w^k\rangle.
\end{aligned}
\end{equation}

One can verify that $W^k({\bf 1}\otimes \bar{x}^k)=0$ holds. Hence, we can add $-W^k({\bf 1}\otimes \bar{x}^k)$ to the right hand side of (\ref{eq:consensus1}) to obtain
\begin{equation}\label{eq:consensus11}
\begin{aligned}
\|x^{k+1}\hspace{-0.05cm} -\hspace{-0.05cm} {\bf 1}\otimes \bar{x}^k\|^2 &\leq
 \| W^k (x^k\hspace{-0.05cm} -\hspace{-0.05cm} {\bf 1}\otimes \bar{x}^k)+\gamma^kr^k \|^2+\|\chi^k\zeta_w^k\|^2\\
&\quad+2\langle W^k  x^k +\gamma^kr^k,  \chi^k\zeta_w^k\rangle.
\end{aligned}
\end{equation}

Now we characterize the first term on the right hand side of (\ref{eq:consensus11}). Using Lemma \ref{Le:rho_2} {on page 5}, we have that there always exists a $T\geq  0$ such that the following inequality holds  for all $k\ge T$:
\[
\begin{aligned}
&\| W^k (x^k\hspace{-0.05cm} -\hspace{-0.05cm} {\bf 1}\otimes \bar{x}^k)+\gamma^kr^k \|^2\\
&\leq \left(\| W^k\|\|x^k-{\bf 1}\otimes \bar{x}^k\|+\gamma^k\|r^k \|\right)^2\\
& \leq \left((1-\chi^k|\rho_2|)\|x^k-{\bf 1}\otimes \bar{x}^k\| +\gamma^k\|r^k \| \right)^2.
\end{aligned}
\]
Next, applying  the inequality
$(a+b)^2\le (1+\epsilon) a^2 + (1+\epsilon^{-1})b^2$, which is valid for any scalars $a,b,$ and $\epsilon>0$, to  the preceding inequality yields
\[
\begin{aligned}
&\| W^k (x^k\hspace{-0.05cm} -\hspace{-0.05cm} {\bf 1}\otimes \bar{x}^k)+\gamma^kr^k \|^2\\
&\leq  (1\hspace{-0.05cm}-\hspace{-0.05cm}\chi^k|\rho_2|)^2(1\hspace{-0.05cm}+\hspace{-0.05cm}\epsilon)\|x^k\hspace{-0.08cm}-\hspace{-0.08cm}{\bf 1}\otimes \bar{x}^k\|^2   +(1\hspace{-0.05cm}+\hspace{-0.05cm}\epsilon^{-1})(\gamma^k)^2\|r^k \|^2
 \end{aligned}
\]
for all $k\ge T$.

Setting $\epsilon=\frac{\chi^k|\rho_2|}{1-\chi^k|\rho_2|}$ (which leads to $(1+\epsilon)=\frac{1}{1-\chi^k|\rho_2|}$ and $1+\epsilon^{-1}=\frac{1}{\chi^k|\rho_2|}$) gives
\begin{equation}\label{eq:consensus2}
\begin{aligned}
&\| W^k (x^k\hspace{-0.05cm} -\hspace{-0.05cm} {\bf 1}\otimes \bar{x}^k)+\gamma^kr^k \|^2\\
&\leq  (1-\chi^k|\rho_2|)\|x^k-{\bf 1}\otimes \bar{x}^k\|^2  +\frac{(\gamma^k)^2}{\chi^k|\rho_2|}\|r^k \|^2.
 \end{aligned}
\end{equation}
Plugging (\ref{eq:consensus2}) into (\ref{eq:consensus11}) yields  for all $k\ge T$:
\begin{equation}\label{eq:consensus3}
\begin{aligned}
&\|x^{k+1}-{\bf 1}\otimes \bar{x}^k\|^2  \leq (1-\chi^k|\rho_2|)\|x^k-{\bf 1}\otimes \bar{x}^k\|^2\\
&\:\:+\frac{(\gamma^k)^2}{\chi^k|\rho_2|}\|r^k \|^2+(\chi^k)^2\|\zeta_w^k\|^2+2\langle W^k x^k+\gamma^kr^k, \chi^k\zeta_w^k\rangle.
\end{aligned}
\end{equation}
Note that when the points $x_i^{k+1}$ are constrained to lie in a convex set $X$, their average $\bar{x}^{k+1}$ also lies in the set $X$, implying
  that the average $\bar{x}^{k+1}$ is the solution to the constrained problem of
minimizing $\sum_{i=1}^{m}\|x_i^{k+1}-x\|^2$ over $x\in X$.
Using this  and the relation $\sum_{i=1}^{m}\|x_i^{k+1}-x\|^2=\|x^{k+1}-{\bf 1}\otimes x\|^2,$
we obtain that
$\|x^{k+1}-{\bf 1}\otimes \bar{x}^{k+1}\|^2\leq \|x^{k+1}-{\bf 1}\otimes \bar{x}^k\|^2$, implying that
for all $k\ge T$,
\begin{equation}\label{eq:consensus4}
\begin{aligned}
&\|x^{k+1}-{\bf 1}\otimes \bar{x}^{k+1}\|^2 \leq (1-\chi^k|\rho_2|)\|x^k-{\bf 1}\otimes \bar{x}^k\|^2 \\
& +\frac{(\gamma^k)^2}{\chi^k|\rho_2|}\|r^k \|^2 +(\chi^k)^2\|\zeta_w^k\|^2+2\langle W^k x^k+\gamma^kr^k, \chi^k\zeta_w^k\rangle.
\end{aligned}
\end{equation}

 Using
$\|\zeta_w^k\|^2\le m \sum_{i=1}^m\sum_{j\neq i} w_{ij}^2\|\zeta_{ij}^k\|^2$, and
the assumption that the DP noise $\zeta_i^k$ has zero mean and variance
$(\sigma_{i}^k)^2$, conditional  on $x_i^k$
(see Assumption~\ref{ass:dp-noise}), taking the conditional expectation, given $\mathcal{F}^k=\{x^0,\,\ldots,x^k\}$,
 we obtain  the following inequality  for all  $k\ge T$ ({note that $ \zeta_w^k$ affects $x^{k+1}$ but it is independent of $x^k$}):
 \begin{equation}\label{eq:consensus4}
\begin{aligned}
&\mathbb{E}\left[\|x^{k+1}\hspace{-0.1cm}-\hspace{-0.1cm}{\bf 1}\otimes \bar{x}^{k+1} \|^2 \mid \mathcal{F}^k\right]\leq (1\hspace{-0.07cm}-\hspace{-0.07cm}\chi^k|\rho_2|)\|x^k\hspace{-0.1cm}-\hspace{-0.1cm}{\bf 1}\otimes \bar{x}^k\|^2 \\
& +\frac{(\gamma^k)^2}{\chi^k|\rho_2|}\|r^k \|^2+ m  (\chi^k)^2\sum_{i=1}^m\sum_{j\neq i} w_{ij}^2(\sigma_j^k)^2.
\end{aligned}
\end{equation}

Therefore, under Assumption \ref{ass:dp-noise} and the conditions for $\chi^k$ and $\gamma^k$ in (\ref{eq:condtions_chi}), we have that the sequence $\|x^k- {\bf 1}\otimes\bar{x}^k\|^2=\sum_{i=1}^{m}\|x_i^k- \bar{x}^k\|^2$ satisfies the conditions for $\{v^k\}$ in Lemma~\ref{Lemma-polyak} {on page 4} (with the index shift),
and hence, converges to zero almost surely. (Note that since the results of Lemma~\ref{Lemma-polyak} are asymptotic, they remain
valid when the starting index is shifted from $k=0$ to $k=T$, for an arbitrary $T\geq0$.)
 Hence, $\lim_{k\to\infty}\|x_i^k - \bar{x}^k\|=0$ {\it a.s.} for all $i\in[m]$.

Moreover, Lemma \ref{Lemma-polyak} also implies  the following relation  {\it a.s.}:
\begin{equation}\label{eq:square}
 \textstyle\sum_{k=0}^{\infty}\chi^k \|x^k- {\bf 1}\otimes \bar{x}^k\|^2 <\infty.
\end{equation}
We invoke the Cauchy–Schwarz inequality to prove the last statement:
\[
\begin{aligned}
&\textstyle\sum_{k=0}^{\infty} \sqrt{\chi^k \|x^k- {\bf 1}\otimes \bar{x}^k\|^2}\sqrt{\frac{(\gamma^k)^2}{\chi^k}}\\
&\qquad\qquad\leq \sqrt{\textstyle\sum_{k=0}^{\infty} \chi^k \|x^k- {\bf 1}\otimes \bar{x}^k\|^2}\sqrt{\textstyle\sum_{k=0}^{\infty} \frac{(\gamma^k)^2}{\chi^k}}.
\end{aligned}
\]
 Noting that the summand in the left hand side of the preceding inequality is actually $\gamma^k\|x^k- {\bf 1}\otimes\bar{x}^k\|$, and the right hand side of the preceding inequality is less than infinity almost surely due to (\ref{eq:square}) and the assumed condition in~(\ref{eq:condtions_chi}), we have
$\sum_{k=0}^{\infty}\gamma^k \|x^k- {\bf 1}\otimes \bar{x}^k\|<\infty$  almost surely.

Further utilizing the relationship
$\sum_{i=1}^{m}\|x_i^k- \bar{x}^k\|\leq    \sqrt{m \|x^k- {\bf 1}\otimes \bar{x}^k\|^2}$ yields the stated result.
\end{proof}


Next, we prove that Algorithm~1 can ensure   $\epsilon$-DP of individual agents' input variables $r_i^k$, even when the number of iterations tends to infinity.

\subsection{Differential-Privacy Analysis}

To analyze the level of DP protection on individual input signals $r_i^k$, we first note that the adjacency defined in Definition \ref{de:adjacency} can be easily adapted  to the constrained consensus  problem by replacing ${\rm col}(f_i,g_i)$ with $r_i$. Specifically, we can define two adjacent constrained consensus  problems as $\mathcal{P}$ and $\mathcal{P}'$, with the only difference being one input signal (represent it as $r_i$ in $\mathcal{P}$ and $r'_i$ in $\mathcal{P}'$ without loss of generality). In line with the second condition of Definition \ref{de:adjacency},  {signals $r_i$  and $r'_i$ are required to have similar terminal behaviors, i.e., $r_i$  and $r'_i$ should converge to each other gradually. We formalize this condition as   the existence of some $C_r>0$ such that}
\begin{equation}\label{eq:r_convergence_rate}
   \|r^k_i-{r'_i}^k\|_1\leq C_r\chi^k
\end{equation}
  holds for all $k\geq 0$.

{For Algorithm 1, an execution is $\mathcal{A}=\{\vartheta^0,\vartheta^1,\ldots\}$ with $\vartheta^k=x^k$. An observation sequence is $\mathcal{O}=\{o^0,o^1,\ldots\}$ with $o^k=y^k$ (note that $y_j^k=x_j^k+\zeta_j^{k}$, see Algorithm 1 for details). } Similar to  the sensitivity definition of constraint-free distributed optimization   in \cite{huang2015differentially}, we define  the sensitivity of a constrained consensus  algorithm as follows: 
\begin{Definition 1}\label{de:sensitivity}
  At each iteration $k$, for any initial state $\vartheta^0$ and  adjacent constrained consensus  problems  $\mathcal{P}$ and $\mathcal{P'}$,  the sensitivity of Algorithm 1 is {($\mathbb{O}$ denotes the set of all possible observation sequences)}
  \begin{equation}
  \Delta^k\triangleq \sup\limits_{\mathcal{O}\in\mathbb{O}}\left\{\sup\limits_{\vartheta\in\mathcal{R}_{\mathcal{P},\vartheta^0}^{-1}(\mathcal{O}),\:\vartheta'\in\mathcal{R}_{\mathcal{P}',\vartheta^0}^{-1}(\mathcal{O})}\hspace{-0.3cm}\|\vartheta^{k}-\vartheta'^{k}\|_1\right\}.
  \end{equation}
\end{Definition 1}

Then, we have the following lemma:
\begin{Lemma 1}\label{Le:Laplacian}
In Algorithm 1, at each iteration $k$, if each agent's DP-noise vector $\zeta_i^k\in\mathbb{R}^d$  consists of $d$ independent Laplace noises with  parameter $\nu^k$   such that $\sum_{k=1}^{T_0}\frac{\Delta^k}{\nu^k}\leq \bar\epsilon$ {holds for some $\bar\epsilon>0$}, then Algorithm 1 achieves $\epsilon$-differential privacy with the cumulative privacy {level} for iterations   $0\leq k\leq T_0$ less than $\bar\epsilon$.
\end{Lemma 1}
\begin{proof}
The proof of the lemma  follows the same line of reasoning  as that  of Lemma~2 in~\cite{huang2015differentially}.
\end{proof}

\begin{Theorem 1}\label{th:DP_Algorithm1}
Under  the conditions of Theorem \ref{Th:consensus_tracking}, if {$\chi^k=\frac{1}{k^s}$ and $\gamma^k=\frac{1}{k^t}$ with $0.5<s<t\leq1$ and $2t-s>1$, and} all elements of $\zeta_i^k$ are drawn independently from  Laplace distribution ${\rm Lap}(\nu^k)$ with $(\sigma_i^k)^2=2(\nu^k)^2$ satisfying Assumption \ref{ass:dp-noise},  then all agents
will reach consensus almost surely. Moreover,
\begin{enumerate}
\item For any finite number of iterations $T$, Algorithm 1 is  $\epsilon$-differentially private with the cumulative privacy budget bounded by $\epsilon\leq \sum_{k=1}^{T}\frac{C_r\varsigma^k}{\nu^k}$   where {$\varsigma^k\triangleq \sum_{p=1}^{k-1} \Pi_{q=p}^{k-1}(1-\chi^q\bar{w}) \gamma^{p-1}\chi^{p-1} +\gamma^{k-1}\chi^{k-1}$, $\bar{w}\triangleq\min_i\{|w_{ii}|\}$, and $C_r$ is from (\ref{eq:r_convergence_rate})};
\item  The cumulative privacy budget is  finite for $T\rightarrow\infty$  when the sequence  $\{\frac{\gamma^k}{\nu^k}\}$ is summable.
\end{enumerate}
\end{Theorem 1}
\begin{proof}
Since  {$\chi^k=\frac{1}{k^s}$ and $\gamma^k=\frac{1}{k^t}$ satisfy the conditions in (\ref{eq:condtions_chi}) in the statement of Theorem 1 under $0.5<s<t\leq 1$ and  $2t-s>1$,} and  the Laplace noise satisfies Assumption \ref{ass:dp-noise}, the convergence result  follows directly from Theorem 1.

To prove the  statements on privacy, we first analyze the sensitivity of Algorithm~1.
Given two adjacent constrained consensus problems $\mathcal{P}$ and $\mathcal{P'}$, for any given fixed observation $\mathcal{O}$ and initial state $x^0$, the sensitivity depends on $\|x^{k}-x'^{k}\|_1$ according to Definition \ref{de:sensitivity}. Since in $\mathcal{P}$ and $\mathcal{P'}$, there is only one input signal that is different, we  represent this different input signal as the  $i$th one, i.e., $r_i$ in $\mathcal{P}$ and $r'_i$ in $\mathcal{P}'$, without loss of generality.
Since the initial conditions, input signals, and observations of $\mathcal{P}$ and $\mathcal{P'}$  are identical for $j\neq i$, we have $x_j^k={x'_j}^k$ for all $j\neq i$ and $k$. Therefore, $\|x^{k}-x'^{k}\|_1$ is always equal to $\|x_i^{k}-{x'_i}^{k}\|_1$.

 According to   Algorithm~1, we can arrive at
 \[
 \begin{aligned}
 \|x_i^{k+1}\hspace{-0.1cm}-\hspace{-0.08cm}{x'_i}^{k+1}\|_1
 \leq  (1\hspace{-0.08cm}-\hspace{-0.08cm}|w_{ii}|\chi^k)\|x_i^k\hspace{-0.08cm}-\hspace{-0.08cm}{x'_i}^k\|_1  +\gamma^k \|r_i^k\hspace{-0.08cm}-\hspace{-0.08cm}{r_i'}^k\|_1,
 \end{aligned}
 \]
where we used the definition $w_{ii}\triangleq-\sum_{j\in\mathbb{N}_i}w_{ij}$ and  the fact that the observations $x_j^k+\zeta_j^k$ and ${x'_j}^k+{\zeta'_j}^k$ are the same.

 Hence, the sensitivity $\Delta^k$ satisfies
 $
 \Delta^{k+1}\leq (1-|w_{ii}|\chi^k)\Delta^{k}+\gamma^k \|r_i^k-{r'_i}^k\|_1.$
 {Using (\ref{eq:r_convergence_rate}), we have
  \begin{equation}\label{eq:delta1}
 \Delta^{k+1}\leq (1-|w_{ii}|\chi^k)\Delta^{k}+C_r\chi^k\gamma^k,
 \end{equation}
 which, implies the first statement by iteration using Lemma \ref{Le:Laplacian}.

 For the infinity horizon result in the second statement, we exploit Lemma~4 in \cite{chung1954stochastic}. More specially, for   the form of $\chi^k$ and $\gamma^k$ in the theorem statement,    Lemma~4 in \cite{chung1954stochastic} implies that (\ref{eq:delta1}) guarantees  the  existence of some $C$ such that $\Delta^{k}$   satisfies $\Delta^{k}<C\gamma^k$. Using  Lemma \ref{Le:Laplacian}, we can easily obtain $\epsilon\leq \sum_{k=1}^{T}\frac{C\gamma^k}{\nu^k}$. Hence, $\epsilon$ will  be finite even when $T$ tends to infinity
 if $\sum_{k=0}^{\infty}\frac{\gamma^k}{\nu^k}<\infty$ holds.
 }
%
\end{proof}

Note that in consensus  applications such as constraint-free distributed optimization, \cite{huang2015differentially} achieves  $\epsilon$-DP by enforcing the sequence $\{\gamma^k\}$  to be summable (geometrically-decreasing, more specifically), which, however, also makes it impossible to ensure accurate convergence. In our approach, by allowing the sequence $\{\gamma^k\}$  to be non-summable (i.e., $\sum_{k=0}^{\infty}\gamma^k=\infty$), we achieve both accurate convergence and finite cumulative privacy budget, even when the number of iterations goes to infinity.  To our knowledge, this is the first constrained consensus  algorithm that achieves both provable  convergence  and rigorous $\epsilon$-DP, even in the infinite time horizon.

\begin{Remark 1}
  To ensure a bounded  cumulative privacy budget $\epsilon=\sum_{k=1}^{\infty}\frac{C\gamma^k}{\nu^k}$, we make the variance  $\nu^k$ of Laplace noise  increase with time (as we require the sequence $\{\frac{\gamma^k}{\nu^k}\}$ to be summable while the sequence $\{\gamma^k\}$ is non-summable).   However, even with an increasing noise  level, our algorithm judiciously makes the
  decaying speed of $\{\chi^k\}$   outweigh the increasing speed of the noise level  $\nu^k$ to successfully ensure that the  actual noise fed into the system, i.e., $\chi^k{\rm Lap}(\nu^k)$,   decays with time. This is the key reason why  Algorithm 1 can ensure every agent's  accurate convergence even in the presence of  noises with an increasing variance. {Moreover, according to Theorem \ref{Th:consensus_tracking}, the convergence is still ensured if we  scale  $\nu^k$ by any constant value $\frac{1}{\epsilon}>0$  to achieve any desired level of $\epsilon$-DP, as long as  $\nu^k$ (with associated variance $(\sigma_i^k)^2=2(\nu^k)^2$) satisfies Assumption \ref{ass:dp-noise}}. 
\end{Remark 1}

\section{Differentially-private  distributed constrained optimization}\label{se:algorithm_2}
In this section, based on the   primal-dual perturbation technique \cite{kallio1994perturbation,kallio1999large}, we present a fully distributed constrained optimization algorithm that can ensure both accurate convergence and rigorous $\epsilon$-DP.
 We assume that the Slater condition holds for problem~\eqref{eq:formulation}.
Similar to \cite{chang2014distributed,belgioioso2020distributed}, we bound the dual updates in a compact  set.  More specifically, we consider a modified version of (\ref{eq:Lagrange_dual}) as follows:
\begin{equation}\label{eq:saddle_point_problem}
 \max_{\lambda\in\mathcal{D}}\left\{\min_{x_i\in\mathcal{X}_i}\mathcal{L}(x_1,\cdots,x_m,\lambda) \right\},
\end{equation}
where
\begin{equation}\label{eq:saddle_point_prolbem_D}
  \hspace{-0.4cm}\mathcal{D}\hspace{-0.07cm}=\hspace{-0.07cm}\left\{\hspace{-0.07cm}\lambda\geq 0\bigg|\|\lambda\|\hspace{-0.07cm}\leq \hspace{-0.07cm}D_\lambda \hspace{-0.07cm}\triangleq \hspace{-0.07cm}\frac{F(\hat{x})-\tilde{q}}{\min_p\{-\sum_{i=1}^{m}g_{ip}(\hat{x}_i)\}}+1 \right\},
\end{equation}
 with $\hat{x}={\rm col}(\hat{x}_1,\ldots,\hat{x}_m)$ being a Slater point of (\ref{eq:formulation}) and $\tilde{q}=\min_{x_i\in\mathcal{X}_i}\mathcal{L}(x_1,\ldots,x_m,\tilde{\lambda})$ being   the dual function value for some arbitrary $\tilde{\lambda}\geq 0$. It has been proven in \cite{nedic2009approximate} that  any dual  optimal  solution $\lambda^\ast$ of (\ref{eq:Lagrange_dual}) satisfies $\|\lambda^\ast\|\leq \frac{F(\hat{x})-\tilde{q}}{\min_p\{-\sum_{i=1}^{m}g_{ip}(\hat{x}_i)\}}$ (i.e., $\|\lambda^\ast\|\leq D_\lambda-1$) and, thus, lies in $\mathcal{D}$. Hence, we can consider the  modified   Lagrangian dual problem (\ref{eq:saddle_point_problem}) instead of the original Lagrangian dual problem (\ref{eq:Lagrange_dual}). This is crucial to obtain provable convergence of the proposed algorithm.  One can verify that (\ref{eq:Lagrange_dual}) and (\ref{eq:saddle_point_problem}) have the same dual solution and attain the same optimal objective value \cite{chang2014distributed}.  Also,  all saddle points of (\ref{eq:Lagrange_dual}) are  saddle points of (\ref{eq:saddle_point_problem}) and vice versa. In fact, the following proposition ensures that the saddle point  of (\ref{eq:saddle_point_problem}) is also an optimal solution to the original optimization problem (\ref{eq:formulation}):

\begin{Proposition 1}\label{Proposition}
 Let the Slater condition hold  and let
$(x_1^\star,\ldots,x_m^\star, \lambda^\star)\in \cX\times{\cal D}$ be a saddle point of (\ref{eq:saddle_point_problem}). Then,
 $(x_1^\star,\ldots,x_m^\star)$ is an optimal solution to problem~(\ref{eq:formulation}) if and only if the following conditions hold
 for all $x_i\in{\cal X}_i$ and all $i\in[m]$:
  \[(x_i-x_i^\star)^T\hspace{-0.1cm} \left(\hspace{-0.05cm}\nabla f_i^T(x_i^\star)\nabla F\hspace{-0.1cm}\left(\textstyle\sum_{j=1}^m f_j(x_j^\star)\right) +\nabla g_i^T(x_i^\star)\lambda^\star\right)\ge 0,\]
\[\textstyle \sum_{i=1}^{m}g_i(x_i^\star)\leq 0, \quad (\lambda^\star)^T\left(\sum_{i=1}^{m}g_i(x_i^\star)\right)=0.
  \]
\end{Proposition 1}
\begin{proof}
 By Proposition 6.1.1.~in \cite{bno2003}, $x^\star= (x_1^\star,\ldots,x_m^\star)$ is an optimal solution to the primal problem if and only if $x^\star$ is primal feasible  (i.e., $\sum_{i=1}^{m}g_i(x_i^\star)\le 0$), and
\[{\cal L}(x^\star,\lambda^\star)=\min_{x\in{\cal X}}
{\cal L}(x,\lambda^\star),\qquad (\lambda^\star)^T\big(\textstyle\sum_{i=1}^{m}g_i(x_i^\star)\big)=0.\]
The condition for the Lagrangian function decouples across variables $x_i$, so it is equivalent to $x_i^\star$ being a minimizer of ${\cal L}(x_1,\ldots,x_m,\lambda^\star)$ with respect to $x_i\in {\cal X}_i$.
The first-order optimality conditions for a minimizer of ${\cal L}(x_1,\ldots,x_m,\lambda^\star)$ over $x_i\in {\cal X}_i$
give the set of relations for $i\in[m]$.
\end{proof}
{
\begin{Remark 1}
 Approaches for locally estimating $D_\lambda$ and ensuring (26) in a  distributed setting have been proposed in the literature of distributed optimization. For example, \cite{belgioioso2020distributed} achieves this goal by letting individual agents build a superset of $\mathcal{D}$; \cite{zhu2011distributed} lets individual agents use a Slater's vector to locally estimate $D_\lambda$. In fact, \cite{belgioioso2020distributed} shows that in practice, agents do not need to know the exact $\mathcal{D}$, and, instead, can just agree on using some superset of it.
\end{Remark 1}
}
Now we are ready to present our distributed constrained optimization algorithm to problem (\ref{eq:formulation}). Besides letting each agent have a local copy $x_i^k$ of the  primal variable, we also  let each agent have a local copy $\lambda_i^k$ of the dual variable. Furthermore, to ensure that the agents collectively optimize the global objective function using local sharing of information, we  let each agent maintain local estimates  of the average  value $\frac{\sum_{i=1}^{m}f_i(x_i^k)}{m}$ and the average constraint-related value
$\frac{\sum_{i=1}^{m}g_i(x_i^k)}{m}$, which are denoted as $y_i^k$ and $z_i^k$, respectively. The detailed algorithm is summarized in Algorithm~2.

\noindent\rule{0.49\textwidth}{0.5pt}
\noindent\textbf{Algorithm 2: Differentially-private distributed constrained
optimization with guaranteed convergence accuracy}
\noindent\rule{0.49\textwidth}{0.5pt}
\begin{enumerate}[wide, labelwidth=!, labelindent=0pt]
    \item[] Parameters: {Deterministic }sequences $\chi^k>0$, $\gamma^k>0$, and $\theta^k>0$.
    \item[] Every agent $i$ maintains a decision variable  $x_i^k$ and a dual variable $\lambda_i^k$, which are initialized randomly in $\mathcal{X}_i$ and $\mathcal{D}$, respectively. Every agent also maintains two auxiliary variables $y_i^k$ and $z_i^k$, which are   initialized as $y_i^0=f_i(x_i^0)$ and $z_i^0=g_i(x_i^0)$, respectively.

    \item[] {\bf for  $k=0,1,\ldots$ do}
    \begin{enumerate}
        \item Every agent generates perturbed primal and dual variables $\alpha_i^k$ and $\beta_i^k$ as follows:
        \begin{equation}\label{eq:perburbed}
         \begin{aligned}
         \alpha_i^{k+1}&\hspace{-0.06cm}=\Pi_{\mathcal{X}_i}\hspace{-0.1cm}\left[x_i^k\hspace{-0.05cm}-\hspace{-0.05cm}\rho_1\left(\nabla f_i^T(x_i^k) \nabla F (my_i^k)\hspace{-0.05cm}+\hspace{-0.05cm}\nabla g_i^T(x_i^k)\lambda_i^k\right)\right],\\
         \beta_i^{k+1}&\hspace{-0.06cm}=\Pi_{\mathcal{D}}\hspace{-0.06cm}\left[\lambda_i^k+\rho_2mz_i^k\right].
         \end{aligned}
        \end{equation}
        \item Every agent  $j$  adds   DP noises   $\zeta_j^{k}$, $\xi_j^k$, and $\vartheta_j^k$ to $\lambda_j^k$, $y_j^k$, and $z_j^k$, respectively, and then  sends the obscured values ${\tilde\lambda_j^k\triangleq}\lambda_j^k+\zeta_j^{k}$, ${\tilde{y}_j^k\triangleq}y_j^k+\xi_j^k$, and ${\tilde{z}_j^k\triangleq}z_j^k+\upsilon_j^k$ to agent
        $i\in\mathbb{N}_j$.
        \item After receiving  ${\tilde{\lambda}_j^k}$, ${\tilde{y}_j^k}$, and ${\tilde{z}_j^k}$  from all $j\in\mathbb{N}_i$, agent $i$ updates its primal and dual variables:
        \begin{equation}\label{eq:primal_dual_update}
        \begin{aligned}
         \hspace{-0.3cm} x_i^{k+1}&\hspace{-0.1cm}=\Pi_{\mathcal{X}_i}\hspace{-0.09cm}\left[x_i^k\hspace{-0.05cm}-\hspace{-0.05cm}\gamma^k\hspace{-0.05cm}\left(\nabla f_i^T(x_i^k) \nabla F (my_i^k)\hspace{-0.05cm}+\hspace{-0.05cm}\nabla g_i^T(x_i^k)\beta_i^{k+1}\right)\right],\\
         \hspace{-0.3cm} \lambda_i^{k+1}&\hspace{-0.1cm}=\Pi_{\mathcal{D}}\hspace{-0.05cm}\left[\lambda_i^k\hspace{-0.05cm}+\hspace{-0.05cm}\chi^k{\textstyle\sum_{j\in \mathbb{N}_i}}\hspace{-0.05cm} w_{ij}({\tilde{\lambda}_j^k}\hspace{-0.05cm}-\hspace{-0.05cm}\lambda_i^k) +\gamma^k g_i(\alpha_i^{k+1})\right],
         \end{aligned}
        \end{equation}
        and auxiliary variables { $y_i^k$ (used in the update of $\alpha_i^k$ and $x_i^k$)  and $z_i^k$ (used in the update of $\beta_i^k$)}:
        \begin{equation}\label{eq:auxiliary_update}
        \begin{aligned}
         y_i^{k+1}&= (1-\theta^k)y_i^k+\chi^k{\textstyle\sum_{j\in \mathbb{N}_i}} w_{ij}({\tilde{y}_j^k}-y_i^k)\\
                  &\qquad   +f_i(x_i^{k+1})-(1-\theta^k)f_i(x_i^k),\\
         z_i^{k+1}&= (1-\theta^k)z_i^k+\chi^k{\textstyle\sum_{j\in \mathbb{N}_i}} w_{ij}({\tilde{z}_j^k}-z_i^k)\\
                  &\qquad   +g_i(x_i^{k+1})-(1-\theta^k)g_i^k(x_i^k).
         \end{aligned}
        \end{equation}
                \item {\bf end}
    \end{enumerate}
\end{enumerate}
\vspace{-0.1cm} \rule{0.49\textwidth}{0.5pt}


 {Compared with the distributed optimization algorithm under inequality constraints without considering privacy in \cite{chang2014distributed}, the proposed algorithm has the following major differences: 1) It introduces  the    sequence $\{\chi^k\}$, which diminishes with time and gradually suppresses the influence of persistent DP noises $\zeta_i^k$, $\xi_i^k$, and $\upsilon_i^k$   on the convergence point of the iterates; 2) It introduces an extra stepsize $\{\theta^k\}$ in the dynamics of the auxiliary variables $y_i^k$ and $z_i^k$, which is necessary to suppress the influence of  DP noises on the local estimation of the average mapping  $\frac{\sum_{i=1}^{m}f_i(x_i^k)}{m}$ and the average constraint-related value
$\frac{\sum_{i=1}^{m}g_i(x_i^k)}{m}$.}
The  sequences  $ \{\gamma^k\}$, $\{\theta^k\}$, and $\{\chi^k\}$
have to be designed appropriately to guarantee the accurate convergence of the iterate vector $x^k\triangleq{\rm col}(x_1^k,\ldots,x_m^k)$ to an optimal solution $x^{\ast}\triangleq {\rm col}(x_1^\ast,\ldots,x_m^\ast)$. {They   are hard-coded into all agents' programs and
need no adjustment/coordination in implementation.}
The persistent DP noise sequences $\{\zeta_i^k\},\, \{\xi_i^k\},\, \{\upsilon_i^k\},\, i\in[m]$ have zero-means and
$\chi^k$-bounded  variances,
 as specified  next.

\begin{Assumption 1}\label{ass:dp-noises-game}
In Algorithm~2, for every $i\in[m]$, the noise sequences $\{\zeta_i^k\}$,
$\{\xi_i^k\}$, and $\{\upsilon_i^k\}$ are zero-mean independent random variables, and independent of $\{x_i^0,\lambda_i^0,z_i^0,y_i^0;i\in[m]\}$.
Also, for every $k$, the noise collection $\{\zeta_j^k,\xi_j^k,\upsilon_j^k; j\in[m]\}$ is independent. The noise variances
$(\sigma_{\zeta,i}^k)^2=\mathbb{E}\left[\|\zeta_i^k\|^2\right]$, $(\sigma_{\xi,i}^k)^2=\mathbb{E}\left[\|\xi_i^k\|^2\right]$, and $(\sigma_{\upsilon,i}^k)^2=\mathbb{E}\left[\|\upsilon_i^k\|^2\right]$
satisfy
\vspace{-0.1cm}
\begin{equation}\label{eq:condition_assumption5}
\begin{aligned}
\hspace{-0.27cm}\sum_{k=0}^\infty(\chi^k)^2\hspace{-0.07cm} \max_{i\in[m]}\hspace{-0.07cm}\left\{\max\hspace{-0.07cm}\left\{(\sigma_{\zeta,i}^k)^2,(\sigma_{\xi,i}^k)^2,(\sigma_{\upsilon,i}^k)^2\right\}\right\}<\infty.
\end{aligned}
\vspace{-0.05cm}
\end{equation}
The initial random vectors $x_i^0$ and $\lambda_i^0$ are such that for all $i$,
$\max\{ \mathbb{E}\left[\|x_i^0\|^2\right],\,
\mathbb{E}\left[\|\lambda_i^0\|^2\right],\, \mathbb{E}\left[\|z_i^0\|^2\right],\,\mathbb{E}\left[\|y_i^0\|^2\right]\}<\infty$.
\end{Assumption 1}

In the following two sections, we prove that Algorithm~2 can 1)~ensure almost sure convergence of all agents to an optimal solution, even in the presence of persistent DP noises, and 2)~achieve rigorous $\epsilon$-DP for individual agents' both cost and constraint functions, even when the number of iterations tends to infinity. To our knowledge, this is the first time that both $\epsilon$-DP and accurate convergence   are achieved in distributed optimization subject to a shared inequality constraint.

\section{Convergence to an Optimal Solution}\label{se:convergence_algo2}
In this section, we will prove that  Algorithm~2 can ensure almost sure  convergence of all agents' iterates to an optimal solution. {Note that different from~\cite{chang2014distributed} which constructs optimal solutions using  running averages of its iterates, we will prove that our iterates $x^k={\rm col}(x_1^k,\ldots,x_m^k)$ will directly converge to an optimal solution to~(\ref{eq:formulation}), which avoids any additional computational and storage overheads  for maintaining additional running averages of iterates.}

\begin{Lemma 1}\label{le:summable}
  Under Assumptions \ref{ass:compact},  \ref{ass:network_function_Lipschitz},  \ref{as:L}, and  \ref{ass:dp-noises-game},  if
   the nonnegative sequences $\{\gamma^k\}$, $\{\chi^k\}$, and $\{\theta^k\}$ in Algorithm 2   satisfy
    \begin{equation}\label{eq:condtions_chi_Algo2}
    \begin{aligned}
     &\quad \sum_{k=0}^{\infty}\chi^k = \infty,\: \sum_{k=0}^{\infty}\theta^k =\infty,\: \sum_{k=0}^{\infty}\gamma^k =\infty,\\  &\sum_{k=0}^{\infty}(\chi^k)^2\hspace{-0.05cm}<\hspace{-0.05cm}\infty,  \sum_{k=0}^{\infty}\frac{(\gamma^k)^2}{\theta^k}\hspace{-0.05cm}<\hspace{-0.05cm}\infty, \sum_{k=0}^{\infty}\frac{(\theta^k)^2}{\chi^k}\hspace{-0.05cm}<\hspace{-0.05cm}\infty,
     \end{aligned}
    \end{equation}
      then,
     the following relations  hold almost surely for $\forall i\in[m]$:
    \begin{equation}\label{eq:lambda}
      {\textstyle\sum_{k=0}^{\infty}}\gamma^k\|\lambda_i^k-\bar{\lambda}^k\|<\infty, \quad \lim_{k\rightarrow \infty}\|\lambda_i^k-\bar{\lambda}^k\|=0,
    \end{equation}
        \begin{equation}\label{eq:y}
       {\textstyle\sum_{k=0}^{\infty}}\gamma^k\|y_i^k-\bar{y}^k\|<\infty, \quad \lim_{k\rightarrow \infty}\|y_i^k-\bar{y}^k\|=0,
    \end{equation}
            \begin{equation}\label{eq:z}
       {\textstyle\sum_{k=0}^{\infty}}\gamma^k\|z_i^k-\bar{z}^k\|<\infty, \quad \lim_{k\rightarrow \infty}\|z_i^k-\bar{z}^k\|=0.
    \end{equation}
\end{Lemma 1}
\begin{proof}
 The relation in (\ref{eq:lambda}) can be obtained by applying Theorem \ref{Th:consensus_tracking} to the dynamics of $\lambda_i^k$ in (\ref{eq:primal_dual_update}). The relations for $y_i^k$ and $z_i^k$ can be obtained by applying the results of Theorem 1 in \cite{wang2022differentially_consensus} to the dynamics of $y_i^k$ and $z_i^k$ in (\ref{eq:auxiliary_update}) (noting that the sequence $\{\gamma^k\}$ decreases faster than $\{\theta^k\}$).
\end{proof}

Next we prove that the perturbation points  $\alpha_i^k$ and $\beta_i^k$ in (\ref{eq:perburbed}) for all agents will converge to their respective ``centralized" counterparts, which are defined below:
\begin{equation}\label{eq:centralized_perturb}
 \begin{aligned}
 \hat\alpha_i^{k+1}\hspace{-0.1cm}&=\hspace{-0.08cm}\Pi_{\mathcal{X}_i}\hspace{-0.08cm}\left[x_i^k\hspace{-0.08cm}-\hspace{-0.08cm}\rho_1\hspace{-0.08cm}\left(\hspace{-0.02cm}\nabla \hspace{-0.04cm} f_i^T\hspace{-0.08cm}(\hspace{-0.03cm}x_i^k\hspace{-0.03cm})\nabla\hspace{-0.08cm} F ({\textstyle\sum_{i=1}^{m}}f_i(\hspace{-0.03cm}x_i^k\hspace{-0.03cm})) \hspace{-0.08cm}+\hspace{-0.08cm}\nabla \hspace{-0.04cm} g_i^T\hspace{-0.08cm}(\hspace{-0.03cm}x_i^k\hspace{-0.03cm})\bar\lambda^k\hspace{-0.04cm}\right)\right],\\
         \hat\beta^{k+1}\hspace{-0.08cm}&=\hspace{-0.08cm}\Pi_{\mathcal{D}}\left[\bar\lambda^k+\rho_2{\textstyle \sum_{i=1}^{m}}g_i(x_i^k)\right],
 \end{aligned}
\end{equation}
where  $\bar{\lambda}^k$ is the average value of all $\lambda_i^k$. We call $\hat\alpha_i^k$ and $\hat\beta^{k}$ the ``centralized" counterparts of $\alpha_i^k$ and $\beta_i^k$ in (\ref{eq:perburbed})  because they are updated based on $\sum_{i=1}^{m}f_i(x_i^k) $, $\sum_{i=1}^{m}g_i(x_i^k)$, and $\bar{\lambda}^k$, which are not directly available to individual agents.

\begin{Lemma 1}\label{Le:alpha_beta}
Under the conditions of Lemma \ref{le:summable}, we have  that  the difference between individual agents' perturbation points $\alpha_i^k$, $\beta_i^k$ in (\ref{eq:perburbed}) and their respective centralized counterparts $\hat{\alpha}^k_i$
and $\hat{\beta}^k$ in~(\ref{eq:centralized_perturb}) converge to 0 almost surely,
i.e.,
\[\lim_{k\to\infty}\|\alpha_i^k-\hat \alpha_i^k\|=0,\qquad
\lim_{k\to\infty}\|\beta_i^k-\hat \beta^k\|{=0}\qquad\hbox{\it a.s.}\]
In addition, the following relations hold almost surely:
  \begin{equation}\label{eq:alpha_summable}
      {\textstyle\sum_{k=0}^{\infty}}\gamma^k\|\alpha_i^k-\hat{\alpha}_i^k\|<\infty,\quad {\textstyle\sum_{k=0}^{\infty}}\gamma^k\|\beta_i^k-\hat{\beta}^k\|<\infty.
  \end{equation}
\end{Lemma 1}
\begin{proof}
We can obtain the dynamics of $\bar{y}^k\triangleq \frac{\sum_{i=1}^{m}y_i^k}{m}$ as follows based on (\ref{eq:auxiliary_update}):
\begin{equation}\label{eq:bar_y1}
\bar{y}^{k+1} = (1-\theta^k)\bar{y}^k-\chi^k \bar\xi_{w}^k +\bar{f}^{k+1}-(1-\theta^k)\bar{f}^k,
\end{equation}
where we defined  $
  \bar\xi_w^k\triangleq \frac{\sum_{i=1}^{m}w_{ii}\xi_i^k}{m},\:\bar{f}^k\triangleq \frac{\sum_{i=1}^{m}f_i(x_i^k)}{m}$, and $ \bar{g}^k\triangleq \frac{\sum_{i=1}^{m}g_i(x_i^k)}{m}$ for notational simplicity.

The relationship in (\ref{eq:bar_y1}) can be rewritten as
 \begin{equation}\label{eq:bar_y2}
\bar{y}^{k+1}- \bar{f}^{k+1}= (1-\theta^k)(\bar{y}^k-\bar{f}^k)+\chi^k \bar\xi_{w}^k.
\end{equation}

Then, following a derivation similar to Theorem~1 of~\cite{wang2022differentially_consensus}, we can prove that, almost surely,  $\|\bar{y}^{k+1}- \bar{f}^{k+1}\|\to0$, $\sum_{k=0}^{\infty}\theta^k\|\bar{y}^k-\bar{f}^k\|^2<\infty$, and $\sum_{k=0}^{\infty}\gamma^k\|\bar{y}^k-\bar{f}^k\|<\infty$.

Further, using the nonexpansiveness property of the projection, we  arrive at the following relation from~(\ref{eq:perburbed}) and~(\ref{eq:centralized_perturb}):
  \begin{equation}\label{eq:alpha_1}
   \begin{aligned}
   &\|\alpha_i^{k+1}-\hat{\alpha}_i^{k+1}\|\\
   &=\left\|\Pi_{\mathcal{X}_i}\left[x_i^k-\rho_1\left(\nabla f_i^T(x_i^k)\nabla  {F}(my_i^k)+\nabla g_i^T(x_i^k)\lambda_i^k\right)\right]\right.\\
   &\left.\quad - \Pi_{\mathcal{X}_i}\left[x_i^k-\rho_1\left(\nabla f_i^T(x_i^k)\nabla  {F}(m\bar{f}^k)+\nabla g_i^T(x_i^k)\bar\lambda^k\right)\right] \right\|\\
   &\leq \rho_1\|\nabla f_i^T(x_i^k)\|\|\nabla  {F}(my_i^k)-\nabla  {F}(m\bar{f}^k)\| \\
   &\quad +\rho_1\|\nabla g_i(x_i^k)\|\|\lambda_i^k-\bar\lambda^k\|\\
   &\leq \rho_1L_g\|\lambda_i^{k}-\bar{\lambda}^{k}\|+m\rho_1G_FL_f\|y_i^{k}-\bar{f}^{k}\|\\
    &\leq \rho_1L_g\|\lambda_i^{k}-\bar{\lambda}^{k}\|+m\rho_1G_FL_f(\|y_i^{k}-\bar{y}^{k}\|+\|\bar{y}^{k}-\bar{f}^{k}\|).
   \end{aligned}
  \end{equation}
  Since $\|\lambda_i^k -\bar{\lambda}^k\|$ and $\|y_i^k-\bar{y}^k\|$ converge {\it a.s.} to 0  (see Lemma \ref{le:summable}), and we have proven that  $\|\bar{y}^k- \bar{f}^k\|$ converges {\it a.s.} to 0, (\ref{eq:alpha_1}) implies that
 $\|\alpha_i^k-\hat{\alpha}_i^k\|$   converges {\it a.s.} to $0$ for all $i\in[m]$.
 Combining the proven result $\sum_{k=0}^{\infty}\gamma^k\|\bar{y}^k-\bar{f}^k\|^2<\infty$, (\ref{eq:y}), (\ref{eq:z}), and (\ref{eq:alpha_1}) yields
 ${\textstyle\sum_{k=0}^{\infty}}\gamma^k\|\alpha_i^k-\hat{\alpha}_i^k\|<\infty$ in (\ref{eq:alpha_summable}).

  Using a similar line of argument, we can prove the statements about $\beta_i^k$.
\end{proof}

To prove that Algorithm~2 ensures almost sure convergence to an optimal solution to problem (\ref{eq:formulation}), we will show that even in the presence of DP noises, 1) the primal-dual iterate pairs $( {x}_1^k,\,\ldots,\, {x}_m^k,\bar{\lambda}^k)$ converge to a saddle point of (\ref{eq:saddle_point_problem}); and 2) the saddle point satisfies the primal-dual optimality conditions in Proposition~\ref{Proposition}.

Using Lemmas \ref{le:x} and   \ref{le:mathcal_L} in the Appendix, we can prove that   the primal-dual iterate pairs $(x_1^k,\,\ldots,\, {x}_m^k,\,\bar{\lambda}^k)$ converge to a saddle point of (\ref{eq:saddle_point_problem}).

\begin{Lemma 1}\label{le:convergence_optimal}
  Under  {Assumptions \ref{ass:compact},   \ref{ass:network_function_Lipschitz},   \ref{as:L}, and  \ref{ass:dp-noises-game}},  if $\rho_1$ satisfies $\rho_1\leq \frac{1}{G_J+D_\lambda G_g}$, the sequences $\{\gamma^k\}$,  $\{\theta^k\}$, and $\{\chi^k\}$ satisfy  $\sum_{k=0}^{\infty}\chi^k =\infty, \sum_{k=0}^{\infty}\theta^k =\infty, \sum_{k=0}^{\infty}\gamma^k =\infty, \sum_{k=0}^{\infty}(\chi^k)^2 < \infty,  \sum_{k=0}^{\infty}\frac{(\gamma^k)^2}{\theta^k} < \infty, \sum_{k=0}^{\infty}\frac{(\theta^k)^2}{\chi^k} < \infty$, then Algorithm 2 guarantees that   $x^k={\rm col}(x_1^k,\cdots,x_m^k)$ and $\lambda_1^k,\ldots,\lambda_m^k$ converge to a saddle point $(x^\ast,\lambda^\ast)$ of (\ref{eq:saddle_point_problem}) almost surely. Namely,
  $
    \lim_{k\rightarrow \infty}\|x_i^k- x^\ast_i\|=0$ and $ \lim_{k\rightarrow \infty}\|\lambda_i^k-  {\lambda}^\ast\|=0    $ hold almost surely
   for all $i\in[m]$ where $ {x}^\ast={\rm col}( {x}_1^\ast,\ldots, {x}_m^\ast)$.
\end{Lemma 1}
\begin{proof}
 Let $ {x}^\star={\rm col}( {x}_1^\star,\ldots, {x}_m^\star)$ and $ \lambda^\star={\rm col}(\lambda_1^\star,\ldots, \lambda_m^\star)$ be an arbitrary saddle point of (\ref{eq:saddle_point_problem}). Clearly we have ${x}^\star\in\mathcal{X}$ and $\lambda^\star\in\mathcal{D}$. Plugging ${x}^\star$ and $\lambda^\star$ into the two inequalities in Lemma \ref{le:x} in the Appendix, and summing the inequalities over $i\in[m]$, we can arrive at
  \begin{equation}\label{eq:convergence1}
  \begin{aligned}
    &\|x^{k+1}-x^\star\|^2+{\textstyle\sum_{i=1}^{m}}\mathbb{E}[\|\lambda_i^{k+1}- \lambda^\star\|^2|\mathcal{F}^k]\\
    &\leq  \|x^{k}-x^\star\|^2+{\textstyle\sum_{i=1}^{m}}\|\lambda_i^{k}- \lambda^\star\|^2+ d^k- 2\gamma^k\times\\
    &\left(\mathcal{L}(x^k,\hat{\beta}^{k+1})\hspace{-0.05cm}-\hspace{-0.05cm}\mathcal{L}(x^\star,\hat{\beta}^{k+1})\hspace{-0.05cm}-\hspace{-0.05cm}\mathcal{L}(\hat{\alpha}^{k+1},\bar{\lambda}^k)\hspace{-0.05cm}+\hspace{-0.05cm}\mathcal{L}(\hat{\alpha}^{k+1},\lambda^\star) \right),
    \end{aligned}
  \end{equation}
  with
 $
  d^k \triangleq (c_1+c_6)(\gamma^k)^2+c_7(\chi^k)^2+c_2 \sum_{i=1}^{m}\gamma^k\|y_i^k-\bar{y}^k\|
    +c_3\sum_{i=1}^{m}\gamma^k\|\beta_i^{k+1}-\hat{\beta}^{k+1}\| + c_4  \sum_{i=1}^{m}\gamma^k\|\lambda_i^k-\bar{\lambda}^k\| + c_5\sum_{i=1}^{m}\gamma^k\|\alpha_i^{k+1}-\hat{\alpha}_i^{k+1}\|.
  $
  The coefficients $c_1-c_7$ are constants given in  (\ref{eq:cs}).
  According to the saddle point definition, we have $
  \mathcal{L}(\hat\alpha^{k+1},{\lambda}^\star)\geq \mathcal{L}(x^\star,{\lambda}^\star)$ and  $\mathcal{L}(x^\star,{\beta}^{k+1})\leq \mathcal{L}(x^\star,{\lambda}^\star)$, implying $\mathcal{L}(\hat\alpha^{k+1},{\lambda}^\star)\geq \mathcal{L}(x^\star,{\beta}^{k+1})$. Therefore, we have the following relationship from (\ref{eq:convergence1}):
    \begin{equation}\label{eq:convergence2}
  \begin{aligned}
    &\|x^{k+1}-x^\star\|^2+{\textstyle\sum_{i=1}^{m}}\mathbb{E}[\|\lambda_i^{k+1}- \lambda^\star\|^2|\mathcal{F}^k]\leq  \|x^{k}-x^\star\|^2+\\
    &{\textstyle\sum_{i=1}^{m}}\|\lambda_i^{k}- \lambda^\star\|^2+    d^k- 2\gamma^k\left(\mathcal{L}(x^k,\hat{\beta}^{k+1})-\mathcal{L}(\hat{\alpha}^{k+1},\bar{\lambda}^k) \right).
    \end{aligned}
  \end{equation}

By Lemma \ref{le:mathcal_L} {in the Appendix}, relation $\mathcal{L}(x^k,\hat{\beta}^{k+1})-\mathcal{L}(\hat{\alpha}^{k+1},\bar{\lambda}^k)\geq 0$ holds under the condition  for $\rho_1$. We can also obtain from Lemma \ref{le:summable} and Lemma \ref{Le:alpha_beta} {on page 9} that $d^k$ is summable under the given conditions for $\gamma^k$, $\theta^k$, and $\chi^k$. Hence, $\|x^{k}-x^\star\|^2+\sum_{i=1}^{m}\|\lambda_i^{k}- \lambda^\star\|^2$ satisfies the conditions for $v^k$ in Lemma \ref{Lemma-polyak1} {on page 4}, and will converge almost surely. (Note that since the results of Lemma~\ref{Lemma-polyak1} are asymptotic, they remain
valid when the starting index is shifted from $k=0$ to $k=T$, for an arbitrary $T\geq0$.) Moreover, $\sum_{k=0}^{\infty}\gamma^k\left(\mathcal{L}(x^k,\hat{\beta}^{k+1})-\mathcal{L}(\hat{\alpha}^{k+1},\bar{\lambda}^k) \right)<\infty$ also holds almost surely. Because $\gamma^k$ is non-summable, we have that the limit inferior of $\mathcal{L}(x^k,\hat{\beta}^{k+1})-\mathcal{L}(\hat{\alpha}^{k+1},\bar{\lambda}^k)$ is zero, i.e., there exists a subsequence  $\{\ell_k\}$  such  that  $\mathcal{L}(x^{\ell_k},\hat{\beta}^{\ell_k+1})-\mathcal{L}(\hat{\alpha}^{\ell_{k+1}},\bar{\lambda}^{\ell_k})$ converges to zero  {\it a.s}. This further implies $\lim_{k\rightarrow \infty}\|x^{\ell_k}-\hat{\alpha}^{\ell_k+1}\|=0$ and $\lim_{k\rightarrow \infty}\|\bar\lambda^{\ell_k}-\hat{\beta}^{\ell_k+1}\|=0$   {\it a.s.} by   Lemma~\ref{le:mathcal_L}  in the Appendix. Given that $\{(x^{\ell_k},\bar\lambda^{\ell_k})\}\subset \mathcal{X}\times \mathcal{D}$ is a bounded sequence, it has a limit point  in the set $\mathcal{X}\times \mathcal{D}$ (since $\cX$ and ${\cal D}$ are closed). Without loss of generality, we can assume that $\{(x^{\ell_k},\bar\lambda^{\ell_k})\}$ converges to a point $( x^\ast,\lambda^\ast)\in \mathcal{X}\times \mathcal{D}$, i.e.,
  \begin{equation}\label{eq:convergence3}
    \lim_{k\rightarrow \infty} x^{\ell_k}= x^\ast,\quad  \lim_{k\rightarrow \infty} \bar\lambda^{\ell_k}= \lambda^\ast\qquad\hbox{\it a.s.}
  \end{equation}

The results in (\ref{eq:convergence3}) imply $\lim_{k\rightarrow \infty} \hat\alpha^{\ell_k+1}= x^{\ell_k}= x^\ast$ and $\lim_{k\rightarrow \infty} \hat\beta^{\ell_k+1}= \bar\lambda^{\ell_k}= \lambda^\ast$ by Lemma \ref{le:mathcal_L} in the Appendix. Using (\ref{eq:centralized_perturb}) and the fact that projection is a continuous mapping \cite{bertsekas2015parallel}, we have
\begin{equation}\label{eq:fixed}
 \begin{aligned}
    x_i^\ast & =\hspace{-0.07cm}\Pi_{\mathcal{X}_i}\hspace{-0.07cm}\left[ x_i^\ast \hspace{-0.07cm}-\hspace{-0.07cm}\rho_1\hspace{-0.07cm}\left(\nabla f_i^T\hspace{-0.07cm}( x_i^\ast) \nabla F({\textstyle\sum_{i=1}^{m}}f_i(  x_i^\ast))\hspace{-0.07cm}+\hspace{-0.07cm}\nabla g_i^T\hspace{-0.07cm}(  x_i^\ast)  \lambda^\ast\right)\right]\\
          \lambda^\ast &=\Pi_{\mathcal{D}}\left[ \lambda^\ast+\rho_2m{\textstyle\sum_{i=1}^{m}}g_i(  x_i^\ast)\right],
 \end{aligned}
\end{equation}
i.e., $( x^\ast, \lambda^\ast)$ is a saddle point of the problem (\ref{eq:saddle_point_problem}).
 Note that
$\|x^{\ell_k}-x^\ast\|^2+{\textstyle\sum_{i=1}^{m}}\|\lambda_i^{\ell_k}- \lambda^\ast\|^2\leq  \|x^{\ell_k}-x^\ast\|^2+{\textstyle\sum_{i=1}^{m}}\left(\|\lambda_i^{\ell_k}-\bar\lambda^{\ell_k}\| +\|\bar\lambda^{\ell_k}-\lambda^\ast\|\right)^2.$
Using the proven result that
 $\lambda_i^{\ell_k}-\bar\lambda^{\ell_k}$ converges to $0$  (see Lemma \ref{le:summable} {on page 9}) and that $\lim_{k\rightarrow \infty} x^{\ell_k}=  x^\ast$ and $\lim_{k\rightarrow \infty} \bar\lambda^{\ell_k}= \lambda^\ast$ (see (\ref{eq:convergence3})), we have $\lim_{k\rightarrow \infty}(\|x^{\ell_k}- x^\ast\|^2+\sum_{i=1}^{m}\|\lambda_i^{\ell_k}-  \lambda^\ast\|^2)=0$. Since we have shown that $\{\|x^{ k}-x^\star\|^2+\sum_{i=1}^{m}\|\lambda_i^{ k}- \lambda^\star\|^2\}$ converges almost surely for any saddle point $(x^\star, \lambda^\star)$, we have  $\lim_{k\rightarrow \infty}(\|x^{k}-  x^\ast\|^2+\sum_{i=1}^{m}\|\lambda_i^{k}-  \lambda^\ast\|^2)=0$ {\it a.s.}
\end{proof}

Lemma \ref{le:convergence_optimal} shows that $(x^k,\lambda^k)$ converges {\it a.s.} to a saddle point $(x^\ast,\lambda^\ast)\in \cX\times {\cal D}$, as $k$ tends to infinity, where
$x^*={\rm col}({x}_1^\ast,\ldots,{x}_m^\ast)$  and, for all $i\in[m]$, the vector $x_i^*$
satisfies
\[ x_i^* =\hspace{-0.07cm}\Pi_{\mathcal{X}_i}\hspace{-0.07cm}
\left[x_i^* \hspace{-0.07cm}-\hspace{-0.07cm}\rho_1\hspace{-0.07cm}
\left(\nabla f_i^T\hspace{-0.07cm}(x_i^*) \nabla F ({\textstyle\sum_{i=1}^{m}}f_i(x_i^*))\hspace{-0.07cm}+\hspace{-0.07cm}\nabla g_i^T\hspace{-0.07cm}(x_i^*) \lambda^*\right)\right].\]
 (see relation~\eqref{eq:fixed}).
By Proposition 1.5.8 in~\cite{facchinei2003finite}, the preceding relation is equivalent to
\[\rho_1(x_i-x_i^*)^T
\left(\nabla f_i^T(x_i^*) \nabla F ({\textstyle\sum_{i=1}^{m}}f_i(x_i^*))+\nabla g_i^T\hspace{-0.07cm}(x_i^*) \lambda^*\right)\ge0.\]
Since $\rho_1>0$, it follows that the first stated relation in~Proposition~\ref{Proposition} is satisfied. We next show that $x^*\in\cX$  satisfies the remaining primal-dual optimality conditions of Proposition \ref{Proposition}.
\begin{Lemma 1}\label{le:second}
 Under the assumptions of Lemma~\ref{le:convergence_optimal}, if $\rho_2 \leq \frac{1}{mC_g}$, we  have the following relation almost surely:
  \begin{equation}\label{eq:second_condition}
    \begin{aligned}
       {\textstyle \sum_{i=1}^{m}}g_i( {x}_i^\ast)     =0. 
    \end{aligned}
  \end{equation}
\end{Lemma 1}

\begin{proof}
From the proof of Lemma \ref{le:convergence_optimal},
$\bar{\lambda}^k$ and $\hat\beta^{k+1}$ converge {\it a.s.} to $\lambda^\ast$, and $x_i^k$ converges {\it a.s.} to $x_i^\ast$, where  $({\rm col}(x_1^\ast,\cdots,x_m^\ast),\lambda^\ast)$  is a saddle point of (\ref{eq:saddle_point_problem}).
In addition,  Lemma \ref{le:summable} and Lemma \ref{Le:alpha_beta} {on page 9} ensure that
 $\lambda_i^k-\bar\lambda^k$ and
$\beta_i^k - \hat\beta^k$ converge {\it a.s.} to 0.  Hence, $\lambda_i^k$ and $\beta_i^k$ also  converge {\it a.s.} to  $\lambda^\ast$.

Moreover, Lemma \ref{le:summable}  shows that
$z_i^k- \bar{z}^k$ converges {\it a.s.} to 0, while
$\bar{z}^k$ converges {\it a.s.} to $\frac{\sum_{i=1}^{m}g_i^k(x_i^k)}{m}$ under the dynamics of (\ref{eq:auxiliary_update}) (see   Theorem~1 of \cite{wang2022differentially_consensus}).
Therefore, using the  fact  that projection is a continuous mapping \cite{bertsekas2015parallel}, the following relation follows from the second equality of (\ref{eq:perburbed}):
\begin{equation}\label{eq:primal_optimal}
 \begin{aligned}
          \lambda^{\ast}&=\Pi_{\mathcal{D}}\left[ \lambda^\ast+\rho_2{\textstyle \sum_{i=1}^{m}}g_i(x_i^\ast)\right].
 \end{aligned}
\end{equation}
It has been proven in \cite{nedic2009approximate} that the optimal dual solution $\lambda^\ast$ of (\ref{eq:Lagrange_dual}) satisfies $\|\lambda^\ast\|\leq D_\lambda-1$  (see the discussions below (\ref{eq:saddle_point_prolbem_D})). Hence, under the condition $\rho_2 \leq \frac{1}{mC_g}$ and the fact that $\|g_i\|\leq C_g$, we have
\[
\| \lambda^\ast+\rho_2{\textstyle \sum_{i=1}^{m}}g_i(x_i^\ast) \|\leq \| \lambda^\ast\|+\rho_2\|{\textstyle \sum_{i=1}^{m}}g_i(x_i^\ast) \|\leq D_\lambda,
\]
implying that $\lambda^\ast+\rho_2{\textstyle \sum_{i=1}^{m}}g_i(x_i^\ast)$ lies in $\mathcal{D}$.
Therefore,  (\ref{eq:primal_optimal}) reduces to
$\lambda^{\ast} =  \lambda^\ast+\rho_2{\textstyle \sum_{i=1}^{m}}g_i(x_i^\ast)$,
which means $\rho_2{\textstyle \sum_{i=1}^{m}}g_i(x_i^\ast)=0$ and ${\textstyle \sum_{i=1}^{m}}g_i(x_i^\ast)=0$.
\end{proof}
 In view of Lemma~\ref{le:second} {on page 10} and the discussion preceding the lemma, we have that
the iterates $\{x_i^k\}$ of Algorithm~2 converge to $x_i^*\in \cX_i$ {\it a.s.} for all $i\in[m]$, where $x^*={\rm col}( {x}_1^*,\ldots, {x}_m^*)$ is an optimal solution of problem~\eqref{eq:formulation}.
Lemma \ref{le:convergence_optimal}
and Lemma \ref{le:second} yield  the main results on convergence, as detailed below.
\begin{Theorem 1}\label{th:algorithm_2_convergence}
 Let {Assumptions \ref{ass:compact},   \ref{ass:network_function_Lipschitz},   \ref{as:L}, and \ref{ass:dp-noises-game}} hold. Also, let $\rho_1>0$ and $\rho_2>0$ satisfy
 $\rho_1\leq \frac{1}{G_J+D_\lambda G_g}$ and $\rho_2 \leq \frac{1}{mC_g}$. If  the sequences $\{\gamma^k\}$,  $\{\theta^k\}$, and $\{\chi^k\}$ satisfy  $\sum_{k=0}^{\infty}\chi^k =\infty,  \sum_{k=0}^{\infty}\theta^k =\infty, \sum_{k=0}^{\infty}\gamma^k =\infty, \sum_{k=0}^{\infty}(\chi^k)^2 < \infty,  \sum_{k=0}^{\infty}\frac{(\gamma^k)^2}{\theta^k} < \infty, \sum_{k=0}^{\infty}\frac{(\theta^k)^2}{\chi^k} < \infty$, then Algorithm 2 guarantees that the iterates  ${\rm col}( {x}_1^k,\cdots, {x}_m^k)$  converge to an optimal solution to  problem (\ref{eq:formulation}) almost surely.
\end{Theorem 1}
\begin{proof}
 The theorem follows directly by combining Proposition~\ref{Proposition}, Lemma \ref{le:convergence_optimal},
and Lemma~\ref{le:second} {on page 10}.
\end{proof}

\begin{Remark 1}\label{re:convergence_speed}
 From Lemma \ref{le:convergence_optimal}, we can see that the convergence speed of all $x^k$ to the optimal solution depends on the convexity properties  of $\mathcal{L}$. If $\mathcal{L}$ is strongly convex with respect to the first argument, then we have that the converging speed of $\|x^k-x^\star\|^2$ in Lemma \ref{le:convergence_optimal} to zero is no worse than $\mathcal{O}(\frac{d^k}{\gamma^k})$ according to Lemma 4 in \cite{chung1954stochastic}, where $d^k$ is given below (\ref{eq:convergence1}). Given that $d^k$ is on the order of $ \mathcal{O}(\frac{(\gamma^k)^2}{\chi^k}) $, we have  $\|x^k-x^\star\|^2$ converging to zero with a rate of $\mathcal{O}(\frac{\gamma^k }{\chi^k})$.
\end{Remark 1}

\section{Algorithm 2 Ensures $\epsilon$-Differential Privacy}\label{se:DP_algo2}

{For Algorithm 2, an execution is $\mathcal{A}=\{\vartheta^0,\vartheta^1,\ldots\}$ with $\vartheta^k={\rm col}(\lambda^k,y^k,z^k)$. An observation sequence is $\mathcal{O}=\{o^0,o^1,\ldots\}$ with $o^k={\rm col}(\tilde\lambda^k,\tilde{y}^k,\tilde{z}^k)$ (note that the symbol $\sim$ represents DP-noise obfuscated information, see Algorithm 2 for details). } Similar to Definition \ref{de:sensitivity}, we define the sensitivity of a distributed constrained optimization algorithm to problem (\ref{eq:formulation}) as follows:
\begin{Definition 1}\label{de:sensitivity_game}
  At each iteration $k$, for any initial state $\vartheta^0$ and any adjacent distributed   optimization problems  $\mathcal{P}$ and $\mathcal{P'}$,  the sensitivity of an optimization algorithm is
  \begin{equation}\label{eq:sensitivity_algo2}
  \Delta^k\triangleq \sup\limits_{\mathcal{O}\in\mathbb{O}}\left\{\sup\limits_{\Theta\in\mathcal{R}_{\mathcal{P},\vartheta^0}^{-1}(\mathcal{O}),\:\Theta'\in\mathcal{R}_{\mathcal{P'},\vartheta^0}^{-1}(\mathcal{O})}\hspace{-0.3cm}\|\Theta^{k}-\Theta'^{k}\|_1\right\},
  \end{equation}
  where $\Theta={\rm col}(\lambda^k,y^k,z^k)$ and {$\mathbb{O}$ denotes the set of all possible observation sequences}.
\end{Definition 1}

Similar to Lemma \ref{Le:Laplacian}, we have the following lemma:
\begin{Lemma 1}\label{Le:Laplacian_game}
In Algorithm 2, at each iteration $k$, if each agent adds noise vectors $\zeta_i^k$, $\xi_i^k$, and $\vartheta_i^k$ to its shared messages $\lambda_i^k$, $y_i^k$, and $z_i^k$, respectively, with every noise vector consisting of  independent Laplace scalar noises with  parameter $\nu^k$, such that $\sum_{k=1}^{T_0}\frac{\Delta^k}{\nu^k}\leq \bar\epsilon$, then the iterative distributed Algorithm 2 is $\epsilon$-differentially private with the cumulative privacy {level} for iterations from $k=0$ to $k=T_0$ less than $\bar\epsilon$.
\end{Lemma 1}
\begin{proof}
The lemma can be obtained following the same line of reasoning of Lemma 2 in  \cite{huang2015differentially}.
\end{proof}


{Before giving the main results, we first use Definition 1 and the guaranteed convergence in Theorem 3 to confine the sensitivity.
Note that when the conditions in the statement of Theorem 3 are satisfied, our algorithm ensures convergence of both $\mathcal{P}$ and $\mathcal{P'}$ to their respective optimal solutions, which are the same under the second requirement in Definition \ref{de:adjacency}. This means that  $\|g_i(\alpha_i^{k+1})-g'_i({\alpha'_i}^{k+1})\|_1=0$, $\|f_i(x_i^{k+1})-f'_i({x'_i}^{k+1})\|_1=0$, and $\|g_i(x_i^{k+1})-g'_i({x'_i}^{k+1})\|_1=0$  will hold when $k$ is sufficiently large   (for the iterates  in both $\mathcal{P}$ and $\mathcal{P}'$ to enter the neighborhood of the optimal solution $B_\delta$ in Definition 1, upon which the evolution in $\mathcal{P}$ and $\mathcal{P}'$ will be identical) using the second condition in Definition \ref{de:adjacency}. Furthermore, the ensured convergence also means that $\|g_i(\alpha_i^{k+1})-g'_i({\alpha'_i}^{k+1})\|_1 $, $\|f_i(x_i^{k+1})-f'_i({x'_i}^{k+1})\|_1 $, and $\|g_i(x_i^{k+1})-g'_i({x'_i}^{k+1})\|_1 $ are always  bounded. Hence, there always exists some constant ${C}$ such that the following relations hold {almost surely} for all $k\geq 0$ under the conditions of Theorem 3:
\begin{equation}\label{eq:C_theorem4}
\begin{aligned}
  &\|g_i(\alpha_i^{k+1})-g'_i({\alpha'_i}^{k+1})\|_1<C \chi^k \theta^k ,\\
  &\|f_i(x_i^{k})\hspace{-0.05cm}-\hspace{-0.05cm}f'_i({x'_i}^{k})\|_1\hspace{-0.1cm}<\hspace{-0.1cm} C \chi^{k}\theta^{k},\:  \|g_i(x_i^{k})\hspace{-0.05cm}-\hspace{-0.05cm}g'_i({x'_i}^{k})\|_1\hspace{-0.1cm} <\hspace{-0.1cm}C \chi^k \theta^{k}.
\end{aligned}
\end{equation}
}

\begin{Theorem 1}\label{th:DP_Algorithm2}
 Under the conditions of Theorem \ref{th:algorithm_2_convergence},  if {$\chi^k=\frac{1}{k^s}$, $\theta^k=\frac{1}{k^u}$, and $\gamma^k=\frac{1}{k^t}$ with $0.5<s<u<t\leq 1$, $2u-s>1$, and $2t-u>1$, and} all elements of $\zeta_i^k$, $\xi_i^k$, and $\vartheta_i^k$ are drawn independently from  Laplace distribution ${\rm Lap}(\nu^k)$ with $(\sigma_i^k)^2=2(\nu^k)^2$ satisfying Assumption \ref{ass:dp-noises-game}, then,  all agents will converge almost surely to an optimal solution. Moreover, 
\begin{enumerate}
\item For any finite number of iterations $T$, Algorithm 1 is  $\epsilon$-differentially private with the cumulative privacy budget bounded by $\epsilon\leq \sum_{k=1}^{T}\frac{ C(\varsigma_\lambda^k+\varsigma_y^k+\varsigma_z^k)}{\nu^k}$  {where $\varsigma_\lambda^k\triangleq \sum_{p=1}^{k-1}(\Pi_{q=p}^{k-1}(1-\bar{L}\chi^{q})) \gamma^{p-1}\chi^{p-1}\theta^{p-1}+ \gamma^{k-1}\chi^{k-1}\theta^{k-1}$, $\varsigma_y^k\triangleq  \sum_{p=1}^{k-1}(\Pi_{q=p}^{k-1}(1-\theta^{q}-\bar{L}\chi^{q})(2-\theta^{p-1})\chi^{p-1}\theta^{p-1})+(2-\theta^{k-1})\chi^{k-1}\theta^{k-1}$,  $\varsigma_z^k\triangleq  \sum_{p=1}^{k-1}(\Pi_{q=p}^{k-1}(1-\theta^{q}-\bar{L}\chi^{q})(2-\theta^{p-1})\chi^{p-1}\theta^{p-1})+(2-\theta^{k-1})\chi^{k-1}\theta^{k-1}$,  $\bar{L}\triangleq\min_i\{|L_{ii}|\}$, and $C$ is from (\ref{eq:C_theorem4})};
\item  The cumulative privacy budget is  finite for $T\rightarrow\infty$  when the sequence  $\{\frac{ \theta^k }{\nu^k}\}$ is summable.
\end{enumerate}
\end{Theorem 1}
\begin{proof}
 Since  {$\chi^k=\frac{1}{k^s}$, $\theta=\frac{1}{k^u}$, and $\gamma^k=\frac{1}{k^t}$ satisfy the conditions in  the statement of Theorem 3 under $0.5<s<u<t\leq 1$, $2u-s>1$, and $2t-u>1$, and}  the Laplace noise satisfies Assumption \ref{ass:dp-noises-game}, the convergence result  follows directly from Theorem 3.

According to the definition of sensitivity in Definition \ref{de:sensitivity_game}, we can obtain $\Delta^k=\Delta^k_\lambda+\Delta^k_y+\Delta^k_z$ where $\Delta^k_\lambda$, $\Delta^k_y$, and $\Delta^k_z$ are obtained by replacing $\Theta^k$ in (\ref{eq:sensitivity_algo2}) with $\lambda^k$, $y^k$, and $z^k$, respectively (note that the norm is $L_1$ norm).

Given two adjacent  optimization problems $\mathcal{P}$ and $\mathcal{P'}$,
 we  represent the different mapping/constraint functions as $(f_i,g_i)$ in $\mathcal{P}$ and $(f'_i,g'_i)$ in $\mathcal{P}'$  without loss of generality.

 Here, we only derive the result for $\Delta^k_\lambda$, but $\Delta^k_y$ and $\Delta^k_z$ can be obtained using the same argument.

 Because the initial conditions, mapping/constraint functions, and observations of $\mathcal{P}$ and $\mathcal{P'}$  are identical for $j\neq i$, we have $\lambda_j^k={\lambda'_j}^k$ for all $j\neq i$ and $k$. Therefore, $\|\lambda^{k}-\lambda'^{k}\|_1$ is always equal to $\|\lambda_i^{k}-{\lambda'_i}^{k}\|_1$.

According to   Algorithm  2, we can arrive at
 \[
 \begin{aligned}
 \|\lambda_i^{k+1}-{\lambda'_i}^{k+1}\|_1\leq& (1-|L_{ii}|\chi^k)\|\lambda_i^k-{\lambda'_i}^k\|\\
 &+\gamma^k\|g_i(\alpha_i^{k+1})-g'_i({\alpha'_i}^{k+1})\|_1,
 \end{aligned}
 \]
 where we have   used the   fact that the observations $\lambda_j^k+\zeta_j^k$ and ${\lambda'_j}^k+{\zeta'_j}^k$ are the same.

 Hence,  $\Delta_\lambda^k$ satisfies
 \[
 \Delta_\lambda^{k+1}\leq (1-|L_{ii}|\chi^k)\Delta_\lambda^{k}+\gamma^k \|g_i(\alpha_i^{k+1})-g'_i({\alpha'_i}^{k+1})\|_1.
 \]
{
and further
 \begin{equation}\label{eq:delta_theorem4}
    \Delta_\lambda^{k+1}\leq (1-|L_{ii}|\chi^k)\Delta_\lambda^{k}+C\gamma^k \chi^k\theta^k
 \end{equation}
 according to (\ref{eq:C_theorem4}).
}

Using a similar line of argument, we can obtain $\Delta^k_y$ and $\Delta^k_z$, and hence, arrive at the first privacy statement by iteration.

{For the infinity horizon result in the second privacy statement, we exploit Lemma~4 in \cite{chung1954stochastic}. More specially, for   the form of $\chi^k$, $\theta^k$, and $\gamma^k$ in the theorem statement,    Lemma~4 in \cite{chung1954stochastic} implies that (\ref{eq:delta_theorem4}) guarantees  the  existence of some $C_\lambda$ such that $\Delta_{\lambda}^{k}$ in (\ref{eq:delta_theorem4}) satisfies $\Delta_{\lambda}^{k}<C_\lambda \theta^k$ (note that $\Delta_{\lambda}^{k}<C'_\lambda\gamma^k\theta^k$ for some $C'_\lambda$ implies the existence of some $C_\lambda$ such that $\Delta_{\lambda}^{k}<C_\lambda \theta^k$). Using a similar line of argument, we can obtain that there exist   $C_y$ and $C_z$ such that
$\Delta_y^{k}<C_y \theta^k$ and $\Delta_z^{k}<C_z \theta^k$ hold.

Using  Lemma \ref{Le:Laplacian_game} {on page 11}, we can easily obtain $\epsilon\leq \sum_{k=1}^{T}\frac{C_\lambda \theta^k+C_y \theta^k+C_z \theta^k}{\nu^k}$. Hence, $\epsilon$ will  be finite even when $T$ tends to infinity
 if $\sum_{k=0}^{\infty}\frac{\theta^k}{\nu^k}<\infty$ holds.}
\end{proof}

\begin{Remark 1}
{It is worth noting that   the sensitivity calculation in the Laplace mechanism is only dependent on the probability density function of observations (see  the derivation of Theorem 3.6 in \cite{dwork2014algorithmic} or the derivation in eqn (9) - eqn (12) in \cite{huang2015differentially}), which is not affected by the  probability-zero event of iterates not converging to an optimal solution (it is well known in probability theory that for a random variable with a continuous probability distribution, adding or removing an outcome of probability zero does not affect its  probability density function). Therefore,  the result of almost sure convergence (i.e., convergence with probability one)  in (\ref{eq:C_theorem4}) can lead to the deterministic result in Theorem 4.}
\end{Remark 1}

Note that in the DP design for coupling-constraint free distributed optimization and Nash equilibrium seeking, to avoid the cumulative privacy budget from growing to infinity,  \cite{huang2015differentially} and \cite{ye2021differentially}  employ a summable stepsize  (geometrically-decreasing stepsize, more specifically), which, however, also makes it impossible to ensure accurate convergence. In contrast, our Algorithm 2 allows the stepsize sequence to be non-summable, which enables    both accurate convergence and $\epsilon$-DP with a finite cumulative privacy budget, even in the infinite time horizon.   To our knowledge, this is the first time that both provable convergence and rigorous $\epsilon$-DP are achieved in distributed optimization subject to inequality constraints.

\begin{Remark 1}
  To ensure the  cumulative privacy budget   to be bounded when  $k\rightarrow \infty$,   we allow  the DP-noise  parameter $\nu^k$ to increase  with time (since we require  $\{\frac{ \theta^k }{\nu^k}\}$ to be summable while $\{\theta^k\}$ is non-summable).    However, we judiciously design the decaying factor $\chi^k$ to make its decreasing speed outweigh the increasing speed of the noise level  $\nu^k$, and hence, ensure that the actual noise fed into the algorithm, i.e.,  $\chi^k{\rm Lap}(\nu^k)$,   decreases with time (see  Assumption~\ref{ass:dp-noises-game}). This is why  Algorithm~2 can ensure  accurate convergence  while achieving $\epsilon$-DP. Moreover, from Theorem~\ref{th:algorithm_2_convergence}, {one can see that the convergence will  not be affected if we scale  $\nu^k$ by any constant value $\frac{1}{\epsilon}>0$  to achieve any desired level of $\epsilon$-DP, as long as the DP-noise parameter $\nu^k$  satisfies
 Assumption~\ref{ass:dp-noises-game}}. 
\end{Remark 1}

{
\begin{Remark 1}
  According to Theorem 4, to increase the level of privacy protection (i.e., reduce the privacy budget), we can use a faster increasing $\{\nu^k\}$ for the Laplace noise ${\rm Lap}(\nu^k)$. However, from Assumption \ref{ass:dp-noises-game},   a faster decreasing $\{\chi^k\}$ has to be used when the noise variance increases faster, leading to a slower convergence speed to an optimal solution according to  Remark \ref{re:convergence_speed}. Therefore, a higher level of privacy protection means a lower speed of convergence   to an optimal solution.
\end{Remark 1}
}

\section{Numerical Simulations}\label{se:simulation}
We evaluate the performance of the  proposed distributed constrained optimization algorithm   using the demand side management problem  in smart grid control considered in   \cite{chang2014distributed}. In the problem, we consider a micro grid system involving   $20$ customers. In addition to paying for the market bid, the customers also have to pay additional cost if there is a deviation between the bid purchased in earlier market settlements and the real-time aggregate load of the customers.   The customers locally control their respective deferrable, non-interruptible loads such as electrical vehicles and washing machines to minimize the cost caused by power imbalance. Compared with the traditional centralized control structure, such a distributed control structure enhances robustness to single point of failure \cite{cai2010decentralized}. It also  avoids collecting  real-time power usage  information of customers, which not only infringes on the privacy of customers but also is not easy for a large scale neighborhood.

We represent the total power bids over a time horizon of length $T$ as $p\triangleq {\rm col}(p_1,\ldots,p_T)$, where $p_t$ denotes the power bid at time $1\leq t\leq T$. Similar to \cite{chang2014distributed}, we model the load profile of customer $i$ as $\Psi_i x_i$, where $\Psi_i\in\mathbb{R}^{T\times T}$ is the coefficient matrix composed of load profiles of customer $i$ and $x_i\in\mathbb{R}^T$  denotes the operation scheduling of the appliances of customer $i$. Due to physical conditions and quality-of-service constraints, $x_i$ is constrained in a compact set $\mathcal{X}_i \subset \mathbb{R}^{T}$. The cost-minimization problem can be formulated as the following optimization problem:
\begin{equation}\label{eq:formulation_example}
\begin{aligned}
  \min_{x_i\in\mathcal{X}_i, i\in[m],z\geq 0}&\left(\pi_p\|z\|^2+\pi_s\left\|z-\textstyle\sum_{i=1}^{m}\Psi_ix_i+p\right\|^2\right)\\
  &{\rm s.t.} \textstyle\sum_{i=1}^{m}\Psi_ix_i-p-z\leq 0,
  \end{aligned}
\end{equation}
where the terms  with $\pi_p$ and $\pi_s$ denote  the costs associated with insufficient and excessive power bids, respectively. In the simulation, we generate  a random network graph that satisfies the connectivity condition, i.e., there is a multi-hop path between each pair of customers.

\begin{figure}
\includegraphics[width=0.45\textwidth]{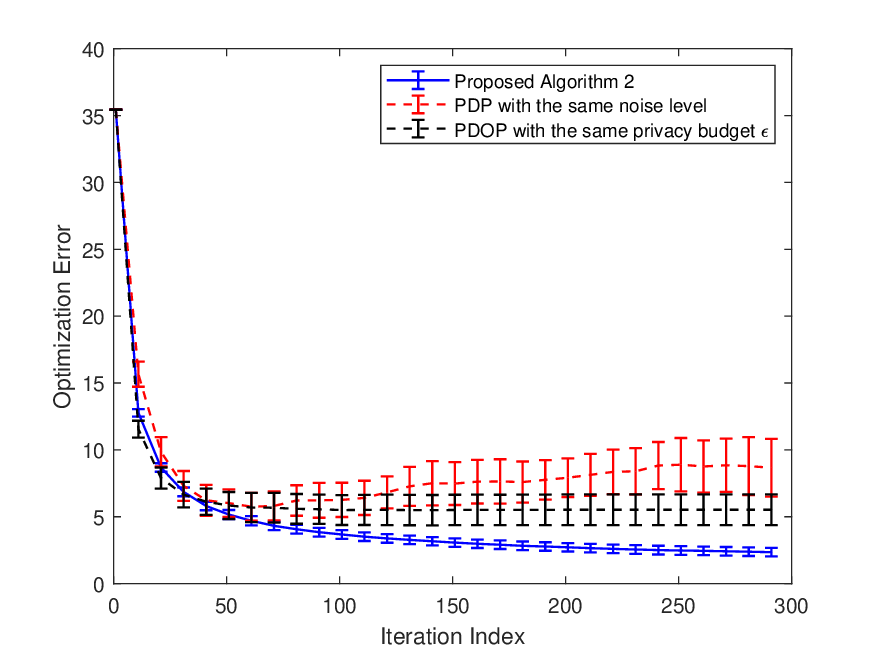}
    \caption{Comparison of Algorithm 2 with the existing distributed constrained optimization algorithm by Chang et al. in \cite{chang2014distributed} (under the same noise level, labeled as PDP) and the differentially private version of the algorithm   by Chang et al. in \cite{chang2014distributed} (using the  DP design in Huang et al. in \cite{huang2015differentially} under the  same cumulative privacy budget $\epsilon$, labeled as PDOP).}
    \label{fig:comparison_algo1}
\end{figure}

To evaluate the performance of the proposed  Algorithm~2, for every agent  (customer)  $i$, we  inject  DP noises $\zeta_i^k$, $\xi_i^k$, and $\upsilon_i^k$ in every shared $\lambda_i^k$, $y_i^k$, and $z_i^k$   in all iterations. Each element of the noise vectors follows Laplace  distribution with parameter $\nu^k=1+0.1k^{0.2}$.  We set the stepsizes and  weakening  sequences as   $\gamma^k=\frac{0.1}{1+0.1k}$, $\theta^k=\frac{0.1}{1+0.1k^{0.96}}$, and $\chi^k =\frac{1}{1+0.1k^{0.9}}$, respectively, which satisfy the conditions in Theorem~3. In the evaluation, we run our algorithm for 100 times, and calculate the average and the variance of  $\|x^k-x^{\ast}\|$  per   iteration index $k$. The result is given by the blue curve and error bars in Fig.~\ref{fig:comparison_algo1}. For comparison, we also run the existing distributed optimization algorithm for a constrained problem  proposed by Chang et al.\ in~\cite{chang2014distributed} under the same  noise level, and the {differentially private version of the algorithm in Chang et al.\ in~\cite{chang2014distributed} under the  DP design in Huang et al.\ in~\cite{huang2015differentially}} using the same cumulative privacy budget $\epsilon$. Note that the DP approach in~\cite{huang2015differentially} addresses  distributed optimization   without shared coupling constraints (it cannot protect the privacy of constraint functions). {We applied its DP mechanism (geometrically decreasing stepsizes and DP noises) to the distributed constrained optimization algorithm in \cite{chang2014distributed} to protect the local cost functions}.   The evolutions of the average error/variance of the  approaches in~\cite{chang2014distributed}  and~\cite{huang2015differentially} are given by  the red and black curves/error bars in Fig.~\ref{fig:comparison_algo1}, respectively. It can be seen that  our proposed Algorithm~2  has a  much better  accuracy.

{

To evaluate the privacy-protection performance of our algorithm, we implemented the attacker models from \cite{melis2019exploiting} and \cite{zhu2019deep} to infer the private information ($\Psi_i$ in (\ref{eq:formulation_example})) in the cost and constraint functions. Note that $\Psi_i$, which is the coefficient matrix composed of load profiles of customer $i$, is contained in the  gradients of both local cost and constraint functions. We implemented different DP noises $\nu^k=\nu_0(1+0.1k^{0.2})$ with $\nu_0$ given by $0,\,0.2,\,0.4,\,0.6,\,0.8$, and $1.0$, respectively. The results {under 300 iterations} are summarized in Table I. It can be seen that with the designed DP protection, the attackers cannot accurately estimate the sensitive information.

\begin{table}
{\caption {Optimization accuracy and attackers' inference errors under different levels of
DP noise}}
 \center
 \vspace{-0.2cm}
{
\begin{tabular}{  |c|c|c|c|c|c|c| }
\hline
  Noise Parameter $\nu_0$ & 0 &  0.2  & 0.4 & 0.6& 0.8&1.0 \\
  \hline
  Optimization Error  &1.75 & 1.84& 1.86& 1.87 & 1.88& 1.88 \\
  \hline
   Inference Error of \cite{melis2019exploiting} &  0.81& 190 & 225 & 247 & 354& 387\\
    \hline
   Inference Error of \cite{zhu2019deep}&  0.12& 104& 132 & 168 & 240 & 259 \\
  \hline
\end{tabular}}
\end{table}
}

\section{Conclusions}\label{se:conclusions}

This paper proposes a differentially-private fully distributed algorithm for distributed optimization subject to a shared inequality constraint. Different from  constraint-free distributed optimization,  the shared inequality constraint  introduces  another  attack surface, and  poses additional challenges to privacy protection. To our knowledge, our approach is the first to achieve differential privacy for both individual cost functions and constraint functions in distributed optimization.  More interestingly, the proposed approach can ensure both accurate convergence to a global optimal solution and rigorous $\epsilon$-differential privacy, which is in sharp contrast  to existing differential privacy solutions for constraint-free distributed optimization that have to sacrifice convergence accuracy for rigorous differential privacy. As a byproduct of the differentially private distributed constrained optimization algorithm, we also propose a new constrained consensus  algorithm that can ensure both provable convergence accuracy and rigorous $\epsilon$-differential privacy, which, to our knowledge, has not been reported before. The convergence analysis also provides a new way to address the entangled dynamics of unbounded DP noises and projection-induced nonlinearity in distributed optimization, which, to our knowledge, has only been addressed separately before.     Numerical simulation results  on a demand response control problem in a smart grid  confirm the effectiveness of the proposed algorithm.

\section*{Appendix}

\begin{Lemma 1}\label{le:x}
  For Algorithm 2, under Assumptions 1-3, there always exists a $T\geq 0$ such that  the following inequalities hold for any $x={\rm col}(x_1,\ldots,x_m)\in\mathcal{X}$ and $\lambda\in\mathcal{D}$ and all $k\geq T$:
  \begin{equation}\label{x^k-x}
  \begin{aligned}
&\sum_{i=1}^{m}\|x_i^{k+1}-x_i\|^2\leq  \sum_{i=1}^{m}\|x_i^{k}-x_i\|^2+c_1(\gamma^k)^2\\
&-2\gamma^k\left(\mathcal{L}(x^k,\hat{\beta}^{k+1})-\mathcal{L}(x,\hat{\beta}^{k+1}) \right)+c_2\sum_{i=1}^{m}\gamma^k\|y_i^k-\bar{y}^k\|\\
    &+c_3\sum_{i=1}^{m}\gamma^k \|\beta_i^{k+1}\hspace{-0.07cm}-\hspace{-0.07cm}\hat\beta^{k+1}\|,
    \end{aligned}
  \end{equation}
    \begin{equation}\label{eq:lambda^k-lambda}
  \begin{aligned}
&\sum_{i=1}^{m}\mathbb{E}\left[\|\lambda_i^{k+1}-\lambda\|^2|\mathcal{F}^k\right]  \leq   \sum_{i=1}^{m}\|\lambda_i^{k} -\hspace{-0.1cm}\lambda\|^2 \\
&+2\gamma^k\left(\mathcal{L}(\hat{\alpha}^{k+1}\hspace{-0.07cm},\bar{\lambda}^k)\hspace{-0.07cm}-\hspace{-0.07cm}\mathcal{L}(\hat{\alpha}^{k+1}\hspace{-0.07cm},\lambda)\right)+ c_4\sum_{i=1}^{m}\gamma^k\|\lambda_i^k-\bar{\lambda}^k\|\\
    &+c_5\sum_{i=1}^{m}\gamma^k\|\alpha_i^{k+1}-\hat{\alpha}_i^{k+1}\|+c_6(\gamma^k)^2+c_7(\chi^k)^2,
    \end{aligned}
  \end{equation}
   where $\mathcal{F}^k=\ \{\lambda_i^\ell,\alpha_i^{\ell+1};\, 0\le \ell\le k\}$ and  $c_1- c_7$ are constants given by
  \begin{equation}\label{eq:cs}
  \begin{aligned}
 & c_1=m(L_fL_F+D_{\lambda} L_g)^2,  c_2=4mD_xL_fG_F,c_3=4D_{\cX}L_g,\\
  &   c_4= 2C_g,c_5=4  D_\lambda  L_g,  c_6= 2 mC_g^2,  c_7=2{\textstyle  \sum_{j\neq i}}w_{ij}^2(\sigma^k_{\zeta,j})^2.
   \end{aligned}
   \end{equation}
\end{Lemma 1}
\begin{proof}
By using the nonexpansiveness property of projection, we can obtain the following relation for any $x={\rm col}(x_1,\ldots, x_m)\in\mathcal{X}$ under Assumptions~\ref{ass:compact} and  \ref{ass:network_function_Lipschitz}:
 \begin{equation}\label{eq:x^k-x0}
 \begin{aligned}
  &\sum_{i=1}^{m}\|x_i^{k+1}-x_i\|^2=\\
  & \sum_{i=1}^{m}\hspace{-0.07cm}\left\|\Pi_{\mathcal{X}_i}\hspace{-0.07cm}\left[x_i^k\hspace{-0.07cm}-\hspace{-0.07cm}\gamma^k\hspace{-0.07cm}\left(\nabla \hspace{-0.07cm} f_i^T\hspace{-0.07cm}(x_i^k)\nabla \hspace{-0.07cm}  {F}(my_i^k)\hspace{-0.07cm}+\hspace{-0.07cm}\nabla g_i^T\hspace{-0.07cm}(x_i^k)\beta_i^{k+1}\right)\right]\hspace{-0.07cm}-\hspace{-0.07cm}x_i\right\|^2\\
  &\leq \sum_{i=1}^{m}\|x_i^k-x_i\|^2+ m(L_fL_F+D_{\lambda} L_g)^2(\gamma^k)^2\\
  &  -2\gamma^k\sum_{i=1}^{m}(x_i^k-x_i)^T \hspace{-0.07cm}\left(\nabla f_i^T(x_i^k)
   \nabla F (my_i^k)\hspace{-0.07cm}+\hspace{-0.07cm}\nabla g_i^T(x_i^k)\beta_i^{k+1}\right).
  \end{aligned}
 \end{equation}
 The last term on the right hand side of the preceding inequality can be bounded as follows:
 \begin{equation}\label{eq:x^k-x}
 \begin{aligned}
 &-2\gamma^k\sum_{i=1}^{m}(x_i^k-x_i)^T \left(\nabla f_i^T(x_i^k) \nabla F (my_i^k)+\nabla g_i^T(x_i^k)\beta_i^{k+1}\right)\\
 &=-2\gamma^k\hspace{-0.1cm}\sum_{i=1}^{m}(x_i^k-x_i)^T \hspace{-0.1cm}\left(\nabla f_i^T(x_i^k) \nabla F (m\bar{y}^k)\hspace{-0.07cm}+\hspace{-0.07cm}\nabla g_i^T(x_i^k)\hat\beta^{k+1}\right)\\
 &-2\gamma^k\sum_{i=1}^{m}(x_i^k-x_i)^T\nabla g_i^T(x_i^k)(\beta_i^{k+1}-\hat\beta^{k+1})\\
 &-2\gamma^k\sum_{i=1}^{m}(x_i^k-x_i)^T \nabla f_i^T(x_i^k) ( \nabla F(my_i^k)-\nabla F (m\bar{y}^k)).
 \end{aligned}
 \end{equation}
 Using the properties  of the gradient of a convex function, we have  for any $x\in\mathcal{X}$,
 \begin{equation}\label{eq:convexity}
   -(x^k-x)^T\mathcal{L}_x(x^k,\hat{\beta}^{k+1})\leq \mathcal{L}(x ,\hat{\beta}^{k+1})-\mathcal{L}(x^k ,\hat{\beta}^{k+1}).
 \end{equation}

Under the Lipschitz and compactness conditions in  Assumptions~\ref{ass:compact} and \ref{ass:network_function_Lipschitz}, combining (\ref{eq:x^k-x}) and (\ref{eq:convexity}) leads to
   \begin{equation}\label{eq:x^k-x2}
 \begin{aligned}
 &-2\gamma^k\sum_{i=1}^{m}(x_i^k-x_i)^T\hspace{-0.1cm} \left(\nabla f_i^T(x_i^k)\nabla  {F}(my_i^k)+\nabla g_i^T(x_i^k)\beta_i^{k+1}\right)\\
 &\leq -2\gamma^k(\mathcal{L}(x^k,\hat{\beta}^{k+1})\hspace{-0.07cm}-\hspace{-0.07cm}\mathcal{L}(x,\hat{\beta}^{k+1}))\hspace{-0.07cm}+\hspace{-0.07cm}c_3\sum_{i=1}^{m}\gamma^k\|\beta_i^{k+1}\hspace{-0.07cm}-\hspace{-0.07cm}\hat\beta^{k+1}\|\\
 &+c_2 \sum_{i=1}^{m}\gamma^k\|y_i^k- \bar{y}^k\|,
 \end{aligned}
 \end{equation}
 with $c_2$ and $c_3$ defined in (\ref{eq:cs}).
Then, by combining (\ref{eq:x^k-x0}) and (\ref{eq:x^k-x2}),   we can arrive at the first  statement  in the lemma.

The second statement in the lemma can be proven following the same line of reasoning. Using  Assumptions~\ref{ass:compact},~\ref{ass:network_function_Lipschitz}, and the update rule of $\lambda_i^k$ in (\ref{eq:primal_dual_update}), we can obtain the following relation
for any $\lambda\in\mathcal{D}$,
\begin{equation}\label{eq:lambda_i-lambda0}
   \begin{aligned}
&\sum_{i=1}^{m}\|\lambda_i^{k+1}-\lambda \|^2   \\
&\leq \sum_{i=1}^{m}\|\sum_{j=1}^{m}\hat w_{ij}\lambda_j^k-\lambda+\chi^k\zeta_{wi}^k+\gamma^kg_i(\alpha_i^{k+1})\|^2\\
&=\sum_{i=1}^{m}(\sum_{j=1}^{m}\|\hat w_{ij}\lambda_j^k-\lambda\|^2+\|\chi^k\zeta_{wi}^k+\gamma^kg_i(\alpha_i^{k+1})\|^2) \\
&\quad +2\sum_{i=1}^{m}(\sum_{j=1}^{m}\hat w_{ij}\lambda_j^k-\lambda)^T (\chi^k\zeta_{wi}^k+\gamma^kg_i(\alpha_i^{k+1})),
   \end{aligned}
 \end{equation}
with $\hat w_{ii}=1-\chi^k\sum_{j=1}^{m}w_{ij}$  and $\hat w_{ij}=\chi^k w_{ij}$ for $j\neq i$, and $\zeta^k_{wi}=\sum_{j=1}^{m}w_{ij}\zeta_j^k$.

Since there always exists some $T\geq 0$ such that $\hat w_{ii}$ is nonnegative, we always have
$
\sum_{i=1}^{m}(\|\sum_{j=1}^{m}\hat w_{ij}\lambda_j^k-\lambda\|^2)\leq   \sum_{j=1}^{m}\|\lambda_j^{k}-\lambda \|^2
$   for all $k\geq T$ (note $\sum_{i=1}^{m}\hat w_{ij}=\sum_{j=1}^{m}\hat w_{ij} =1$ for $i,j\in[m]$).

Therefore, we can arrive at
\[
\begin{aligned}
&\sum_{i=1}^{m}(\sum_{j=1}^{m}\hat w_{ij}\lambda_j^k-\lambda)^T (\gamma^kg_i(\alpha_i^{k+1}))\\
&=\gamma^k\hspace{-0.07cm}\sum_{i=1}^{m}( \bar{\lambda}^k\hspace{-0.07cm}-\hspace{-0.07cm}\lambda)^T\hspace{-0.07cm} g_i(\alpha_i^{k+1})\hspace{-0.07cm}+\hspace{-0.07cm}\gamma^k\sum_{i=1}^{m}(\sum_{j=1}^{m}\hat w_{ij}\lambda_j^k\hspace{-0.07cm}-\hspace{-0.07cm}\bar{\lambda}^k)^T g_i(\alpha_i^{k+1})\\
&\leq \gamma^k\sum_{i=1}^{m}( \bar{\lambda}^k-\lambda)^T g_i(\alpha_i^{k+1})+\gamma^kC_g\sum_{i=1}^{m}\|\lambda_i^k-\bar{\lambda}^k\|.
\end{aligned}
\]
We take the conditional expectation with respect to $\mathcal{F}^k=\ \{\lambda_i^\ell,\alpha_i^{\ell+1};\, 0\le \ell\le k\}$ in (\ref{eq:lambda_i-lambda0}) to obtain
  \begin{equation}\label{eq:lambda_i-lambda}
   \begin{aligned}
&\sum_{i=1}^{m}\mathbb{E}\left[\|\lambda_i^{k+1}-\lambda \|^2|\mathcal{F}^k\right]   \\
    &\leq  \sum_{i=1}^{m} \|\lambda_i^{k} - \lambda \|^2  +2(\chi^k)^2\sum_{i=1}^{m}\sum_{j\neq i}w_{ij}^2(\sigma^k_{\zeta,j})^2+2(\gamma^k)^2 mC_g^2\\
    & +2\gamma^k\hspace{-0.1cm}\sum_{i=1}^{m}(\bar\lambda^k\hspace{-0.07cm}-\hspace{-0.07cm}\lambda)^T g_i(\alpha_i^{k+1})+ 2\gamma^kC_g\sum_{i=1}^{m}\|\lambda_i^k-\bar{\lambda}^k\|.
   \end{aligned}
 \end{equation}

 The second last term on the right hand side of the preceding inequality can be bounded as follows:
 \begin{equation}\label{eq:lambda_i-lambda2}
   \begin{aligned}
     &2\gamma^k\sum_{i=1}^{m}(\bar\lambda^k-\lambda)^T g_i(\alpha_i^{k+1}) =2\gamma^k\sum_{i=1}^{m}(\bar\lambda^k-\lambda)^T g_i(\hat\alpha_i^{k+1})  \\
     &\quad +2\gamma^k\sum_{i=1}^{m}(\bar\lambda^k-\lambda)^T (g_i(\alpha_i^{k+1})-g_i(\hat\alpha_i^{k+1}))\\
     &\leq 2\gamma^k\hspace{-0.07cm}\sum_{i=1}^{m}(\bar\lambda^k\hspace{-0.07cm}-\hspace{-0.07cm}\lambda)^T g_i(\hat\alpha_i^{k+1})  +4D_\lambda L_g\gamma^k\sum_{i=1}^{m}\| \alpha_i^{k+1}\hspace{-0.07cm}-\hspace{-0.07cm} \hat\alpha_i^{k+1}\|.
   \end{aligned}
 \end{equation}
Then, we can arrive at the second statement of the lemma by using   the definition   $\mathcal{L}_\lambda(\hat{\alpha}^{k+1},\bar{\lambda}^k)=\sum_{i=1}^{m}g_i(\hat\alpha_i^{k+1})$  and  the  relationship $ (\bar{\lambda}^k-\lambda)^T\mathcal{L}_\lambda(\hat{\alpha}^{k+1},\bar{\lambda}^k)= \mathcal{L}(\hat{\alpha}^{k+1},\bar{\lambda}^k)-\mathcal{L}(\hat{\alpha}^{k+1}, \lambda )$.
\end{proof}

\begin{Lemma 1}\label{le:mathcal_L}
   For Algorithm 2, under  Assumptions  \ref{ass:compact} and   \ref{ass:network_function_Lipschitz}, the following relationship always holds:
  \begin{equation}\label{eq:perturbation_points}
    \begin{aligned}
      &\mathcal{L}(x^k,\hat{\beta}^{k+1})-\mathcal{L}(\hat{\alpha}^{k+1},\bar\lambda^k)\geq\\
      & (\frac{1}{\rho_1}-G_J-D_\lambda G_g)\|x^k-\hat\alpha^{k+1}\|^2+\frac{1}{\rho_2}\|\bar{\lambda}^k-\hat{\beta}^{k+1}\|^2.
    \end{aligned}
  \end{equation}
\end{Lemma 1}
\begin{proof}
 We first note that (\ref{eq:centralized_perturb}) implies
$\hat\alpha_i^{k+1}=\arg\min_{\alpha_i\in\mathcal{X}_i}\|\alpha_i-x_i^k+\rho_1\mathcal{L}_{x_i}(x^k,\bar\lambda^k)\|^2
  $.
  Hence, for any $x_i\in \mathcal{X}_i$, we have $
  (x_i-\hat\alpha_i^{k+1})^T(\hat\alpha_i^{k+1}-x_i^k+\rho_1\mathcal{L}_{x_i}(x^k,\bar\lambda^k))\geq 0$.
  Setting $x_i$ in the preceding inequality as $x_i^k$ leads to
    $
  \|x_i^k-\hat\alpha_i^{k+1}\|^2\leq  \rho_1(x_i^k-\hat\alpha_i^{k+1})^T \mathcal{L}_{x_i}(x^k,\bar\lambda^k)
  $
  and further
      $
  \|x^k-\hat\alpha^{k+1}\|^2\leq  \rho_1(x^k-\hat\alpha^{k+1})^T \mathcal{L}_{x}(x^k,\bar\lambda^k)
  $
  by summing over $i\in[m]$.

  By replacing $\mathcal{L}_{x}(x^k,\bar\lambda^k)$ in the preceding inequality with $\mathcal{L}_{x}(x^k,\bar\lambda^k)=\mathcal{L}_{x}(x^k,\bar\lambda^k)+\mathcal{L}_{x}(\hat{\alpha}^{k+1},\bar\lambda^k)-\mathcal{L}_{x}(\hat{\alpha}^{k+1},\bar\lambda^k)$, we can arrive at
  \[
  \begin{aligned}
  &\frac{1}{\rho_1}\|x^k-\hat\alpha^{k+1}\|^2\leq  (x^k-\hat\alpha^{k+1})^T \mathcal{L}_{x}(\hat{\alpha}^{k+1},\bar\lambda^k)\\
  &- (x^k-\hat\alpha^{k+1})^T (\mathcal{L}_{x}(\hat{\alpha}^{k+1},\bar\lambda^k)-\mathcal{L}_{x}(x^k,\bar\lambda^k)),
  \end{aligned}
  \]
 implying that
   \[
  \begin{aligned}
  &  (x^k-\hat\alpha^{k+1})^T \mathcal{L}_{x}(\hat{\alpha}^{k+1},\bar\lambda^k) \geq \frac{1}{\rho_1}\|x^k-\hat\alpha^{k+1}\|^2 \\
  &- (x^k-\hat\alpha^{k+1})^T (\mathcal{L}_{x}(x^k,\bar\lambda^k)-\mathcal{L}_{x}(\hat{\alpha}^{k+1},\bar\lambda^k)).
  \end{aligned}
  \]
  The last inequality further implies
  \begin{equation}\label{eq:bound_Lx0}
     \begin{aligned}
  &  (x^k-\hat\alpha^{k+1})^T \mathcal{L}_{x}(\hat{\alpha}^{k+1},\bar\lambda^k) \geq \frac{1}{\rho_1}\|x^k-\hat\alpha^{k+1}\|^2 \\
  &- \|x^k-\hat\alpha^{k+1}\|\|  \mathcal{L}_{x}(x^k,\bar\lambda^k)-\mathcal{L}_{x}(\hat{\alpha}^{k+1},\bar\lambda^k)\|.
  \end{aligned}
  \end{equation}

  Using  Assumptions \ref{ass:compact} and   \ref{ass:network_function_Lipschitz}, and the expression of $\mathcal{L}_x$ in (\ref{eq:L_x_i}), we can bound $\|  \mathcal{L}_{x}(x^k,\bar\lambda^k)-\mathcal{L}_{x}(\hat{\alpha}^{k+1},\bar\lambda^k)\|$ as follows:
  \begin{equation}\label{eq:bound_Lx}
  \begin{aligned}
    &\|\mathcal{L}_{x}(x^k,\bar\lambda^k)-\mathcal{L}_{x}(\hat{\alpha}^{k+1},\bar\lambda^k)\|\leq \|\nabla J(x^k)-\nabla J(\hat{\alpha}^{k+1})\|\\
    &\quad +\|\bar\lambda^k\| \left\|\left[\begin{array}{c} \nabla g_1^T(x_1^k)-\nabla g_1^T(\hat\alpha_1^{k+1})\\ \vdots\\\nabla g_m^T(x_m^k)-\nabla g_m^T(\hat\alpha_m^{k+1}) \end{array} \right]\right\|_F\\
    &\leq (G_J+D_\lambda G_g)\|x^k-\hat\alpha^{k+1}\|,
  \end{aligned}
  \end{equation}
  where $\|\cdot\|_F$ denotes the Frobenius norm.

  Plugging (\ref{eq:bound_Lx}) into (\ref{eq:bound_Lx0}) yields
    \begin{equation}\label{eq:bound_Lx1}
     \begin{aligned}
   (x^k\hspace{-0.06cm}-\hspace{-0.06cm}\hat\alpha^{k+1})^T \hspace{-0.06cm}\mathcal{L}_{x}(\hat{\alpha}^{k+1},\bar\lambda^k)\hspace{-0.06cm} \geq\hspace{-0.06cm} (\frac{1}{\rho_1}\hspace{-0.06cm}-\hspace{-0.06cm} G_J\hspace{-0.06cm}-\hspace{-0.06cm}D_\lambda G_g)\|x^k\hspace{-0.06cm}-\hspace{-0.06cm}\hat\alpha^{k+1}\|^2.
  \end{aligned}
  \end{equation}
  By the convexity of $\mathcal{L}$ with respect to $x$, we have $\mathcal{L}(x^k,\bar{\lambda}^k)-\mathcal{L}(\hat{\alpha}^{k+1},\bar{\lambda}^k)\geq (x^k-\hat\alpha^{k+1})^T \mathcal{L}_{x}(\hat{\alpha}^{k+1},\bar\lambda^k)$, which  implies
      \begin{equation}\label{eq:bound_Lx2}
     \begin{aligned}
  &  \mathcal{L}(x^k,\bar{\lambda}^k)\hspace{-0.07cm}-\hspace{-0.07cm}\mathcal{L}(\hat{\alpha}^{k+1},\bar{\lambda}^k)  \geq (\frac{1}{\rho_1}\hspace{-0.05cm}-\hspace{-0.05cm} G_J\hspace{-0.05cm}-\hspace{-0.05cm}D_\lambda G_g)\|x^k-\hat\alpha^{k+1}\|^2.
  \end{aligned}
  \end{equation}

  Similarly, by the update rule of $\hat\beta^k$ in (\ref{eq:centralized_perturb}), we know $\hat\beta^k=\arg\min_{\beta\in\mathcal{D}}\left\| \beta-\bar\lambda^k-\rho_2\sum_{i=1}^{m}g_i(x_i^k) \right\|^2$. Hence, using the optimality condition and the fact that $\mathcal{L}$ is linear in $\lambda$, we arrive at
  \begin{equation}\label{eq:L_lambda}
    \begin{aligned}
&\mathcal{L}(x^k,\hat{\beta}^{k+1})-\mathcal{L}(x^k,\bar{\lambda}^{k})\\
 &=-(\bar\lambda^k-\hat\beta^{k+1})^T\sum_{i=1}^{m}g_i(x_i^k) \geq \frac{1}{\rho_2}\|\bar\lambda^k-\hat\beta^{k+1}\|^2.
    \end{aligned}
  \end{equation}
  Summing (\ref{eq:bound_Lx2}) and (\ref{eq:L_lambda}) finishes the proof of the lemma.
\end{proof}

{\footnotesize
\bibliography{reference2}}

\begin{thebibliography}{10}
\providecommand{\url}[1]{#1}
\csname url@samestyle\endcsname
\providecommand{\newblock}{\relax}
\providecommand{\bibinfo}[2]{#2}
\providecommand{\BIBentrySTDinterwordspacing}{\spaceskip=0pt\relax}
\providecommand{\BIBentryALTinterwordstretchfactor}{4}
\providecommand{\BIBentryALTinterwordspacing}{\spaceskip=\fontdimen2\font plus
\BIBentryALTinterwordstretchfactor\fontdimen3\font minus
  \fontdimen4\font\relax}
\providecommand{\BIBforeignlanguage}[2]{{%
\expandafter\ifx\csname l@#1\endcsname\relax
\typeout{** WARNING: IEEEtran.bst: No hyphenation pattern has been}%
\typeout{** loaded for the language `#1'. Using the pattern for}%
\typeout{** the default language instead.}%
\else
\language=\csname l@#1\endcsname
\fi
#2}}
\providecommand{\BIBdecl}{\relax}
\BIBdecl

\bibitem{chang2014distributed}
T.-H. Chang, A.~Nedi{\'c}, and A.~Scaglione, ``Distributed constrained
  optimization by consensus-based primal-dual perturbation method,'' \emph{IEEE
  Transactions on Automatic Control}, vol.~59, no.~6, pp. 1524--1538, 2014.

\bibitem{notarnicola2019constraint}
I.~Notarnicola and G.~Notarstefano, ``Constraint-coupled distributed
  optimization: A relaxation and duality approach,'' \emph{IEEE Transactions on
  Control of Network Systems}, vol.~7, no.~1, pp. 483--492, 2019.

\bibitem{patrascu2018convergence}
A.~Patrascu and I.~Necoara, ``On the convergence of inexact projection primal
  first-order methods for convex minimization,'' \emph{IEEE Transactions on
  Automatic Control}, vol.~63, no.~10, pp. 3317--3329, 2018.

\bibitem{hershberger2001distributed}
D.~E. Hershberger and H.~Kargupta, ``Distributed multivariate regression using
  wavelet-based collective data mining,'' \emph{Journal of Parallel and
  Distributed Computing}, vol.~61, no.~3, pp. 372--400, 2001.

\bibitem{notarstefano2019distributed}
G.~Notarstefano, I.~Notarnicola, A.~Camisa \emph{et~al.}, ``Distributed
  optimization for smart cyber-physical networks,'' \emph{Foundations and
  Trends in Systems and Control}, vol.~7, no.~3, pp. 253--383, 2019.

\bibitem{yang2010distributed}
B.~Yang and M.~Johansson, ``Distributed optimization and games: A tutorial
  overview,'' \emph{Networked Control Systems}, pp. 109--148, 2010.

\bibitem{su2022convergence}
Y.~Su, Z.~Wang, M.~Cao, M.~Jia, and F.~Liu, ``Convergence analysis of dual
  decomposition algorithm in distributed optimization: Asynchrony and
  inexactness,'' \emph{IEEE Transactions on Automatic Control}, 2022.

\bibitem{zhu2011distributed}
M.~Zhu and S.~Mart{\'\i}nez, ``On distributed convex optimization under
  inequality and equality constraints,'' \emph{IEEE Transactions on Automatic
  Control}, vol.~57, no.~1, pp. 151--164, 2011.

\bibitem{yuan2011distributed}
D.~Yuan, S.~Xu, and H.~Zhao, ``Distributed primal--dual subgradient method for
  multiagent optimization via consensus algorithms,'' \emph{IEEE Transactions
  on Systems, Man, and Cybernetics, Part B (Cybernetics)}, vol.~41, no.~6, pp.
  1715--1724, 2011.

\bibitem{zhang2019admm}
C.~Zhang, M.~Ahmad, and Y.~Wang, ``{ADMM} based privacy-preserving
  decentralized optimization,'' \emph{IEEE Transactions on Information
  Forensics and Security}, vol.~14, no.~3, pp. 565--580, 2019.

\bibitem{huang2015differentially}
Z.~Huang, S.~Mitra, and N.~Vaidya, ``Differentially private distributed
  optimization,'' in \emph{Proceedings of the 2015 International Conference on
  Distributed Computing and Networking}, 2015, pp. 1--10.

\bibitem{burbano2019inferring}
D.~A. Burbano-L, J.~George, R.~A. Freeman, and K.~M. Lynch, ``Inferring private
  information in wireless sensor networks,'' in \emph{IEEE International
  Conference on Acoustics, Speech and Signal Processing}, 2019, pp. 4310--4314.

\bibitem{wang2022quantization}
Y.~Wang and T.~Ba{\c{s}}ar, ``Quantization enabled privacy protection in
  decentralized stochastic optimization,'' \emph{IEEE Transactions on Automatic
  Control}, 2022.

\bibitem{wang2022tailoring}
Y.~Wang and A.~Nedi{\'c}, ``Tailoring gradient methods for
  differentially-private distributed optimization,'' \emph{IEEE Transactions on
  Automatic Control}, 2023.

\bibitem{wang2022decentralized}
Y.~Wang and H.~V. Poor, ``Decentralized stochastic optimization with inherent
  privacy protection,'' \emph{IEEE Transactions on Automatic Control}, 2022.

\bibitem{falsone2017dual}
A.~Falsone, K.~Margellos, S.~Garatti, and M.~Prandini, ``Dual decomposition for
  multi-agent distributed optimization with coupling constraints,''
  \emph{Automatica}, vol.~84, pp. 149--158, 2017.

\bibitem{tjell2019privacy}
K.~Tjell and R.~Wisniewski, ``Privacy preservation in distributed optimization
  via dual decomposition and admm,'' in \emph{2019 IEEE 58th Conference on
  Decision and Control (CDC)}.\hskip 1em plus 0.5em minus 0.4em\relax IEEE,
  2019, pp. 7203--7208.

\bibitem{han2021privacy}
D.~Han, K.~Liu, H.~Sandberg, S.~Chai, and Y.~Xia, ``Privacy-preserving dual
  averaging with arbitrary initial conditions for distributed optimization,''
  \emph{IEEE Transactions on Automatic Control}, vol.~67, no.~6, pp.
  3172--3179, 2021.

\bibitem{zhang2018enabling}
C.~Zhang and Y.~Wang, ``Enabling privacy-preservation in decentralized
  optimization,'' \emph{IEEE Transactions on Control of Network Systems},
  vol.~6, no.~2, pp. 679--689, 2018.

\bibitem{freris2016distributed}
N.~M. Freris and P.~Patrinos, ``Distributed computing over encrypted data,'' in
  \emph{2016 54th Annual Allerton Conference on Communication, Control, and
  Computing (Allerton)}.\hskip 1em plus 0.5em minus 0.4em\relax IEEE, 2016, pp.
  1116--1122.

\bibitem{lu2018privacy}
Y.~Lu and M.~Zhu, ``Privacy preserving distributed optimization using
  homomorphic encryption,'' \emph{Automatica}, vol.~96, pp. 314--325, 2018.

\bibitem{lou2017privacy}
Y.~Lou, L.~Yu, S.~Wang, and P.~Yi, ``Privacy preservation in distributed
  subgradient optimization algorithms,'' \emph{IEEE Transactions on
  Cybernetics}, vol.~48, no.~7, pp. 2154--2165, 2017.

\bibitem{gade2018private}
S.~Gade and N.~H. Vaidya, ``Private optimization on networks,'' in
  \emph{American Control Conference}.\hskip 1em plus 0.5em minus 0.4em\relax
  IEEE, 2018, pp. 1402--1409.

\bibitem{gao2022dynamics}
H.~Gao, Y.~Wang, and A.~Nedi{\'c}, ``Dynamics based privacy preservation in
  decentralized optimization,'' \emph{Automatica}, vol. 151, p. 110878, 2023.

\bibitem{wang2022decentralized1}
Y.~Wang and A.~Nedi\'c, ``Decentralized gradient methods with time-varying
  uncoordinated stepsizes: Convergence analysis and privacy design,''
  \emph{IEEE Transactions on Automatic Control}, 2022.

\bibitem{dwork2014algorithmic}
C.~Dwork, A.~Roth \emph{et~al.}, ``The algorithmic foundations of differential
  privacy.'' \emph{Foundations and Trends in Theoretical Computer Science},
  vol.~9, no. 3-4, pp. 211--407, 2014.

\bibitem{han2016differentially}
S.~Han, U.~Topcu, and G.~J. Pappas, ``Differentially private distributed
  constrained optimization,'' \emph{IEEE Transactions on Automatic Control},
  vol.~62, no.~1, pp. 50--64, 2016.

\bibitem{hale2017cloud}
M.~T. Hale and M.~Egerstedt, ``Cloud-enabled differentially private multiagent
  optimization with constraints,'' \emph{IEEE Transactions on Control of
  Network Systems}, vol.~5, no.~4, pp. 1693--1706, 2017.

\bibitem{wang2017differential}
Y.~Wang, Z.~Huang, S.~Mitra, and G.~E. Dullerud, ``Differential privacy in
  linear distributed control systems: Entropy minimizing mechanisms and
  performance tradeoffs,'' \emph{IEEE Transactions on Control of Network
  Systems}, vol.~4, no.~1, pp. 118--130, 2017.

\bibitem{zhang2019recycled}
X.~Zhang, M.~M. Khalili, and M.~Liu, ``Recycled {ADMM}: Improving the privacy
  and accuracy of distributed algorithms,'' \emph{IEEE Transactions on
  Information Forensics and Security}, vol.~15, pp. 1723--1734, 2019.

\bibitem{he2020differential}
J.~He, L.~Cai, and X.~Guan, ``Differential private noise adding mechanism and
  its application on consensus algorithm,'' \emph{IEEE Transactions on Signal
  Processing}, vol.~68, pp. 4069--4082, 2020.

\bibitem{cortes2016differential}
J.~Cort{\'e}s, G.~E. Dullerud, S.~Han, J.~Le~Ny, S.~Mitra, and G.~J. Pappas,
  ``Differential privacy in control and network systems,'' in \emph{IEEE 55th
  Conference on Decision and Control (CDC)}, 2016, pp. 4252--4272.

\bibitem{xiong2020privacy}
Y.~Xiong, J.~Xu, K.~You, J.~Liu, and L.~Wu, ``Privacy preserving distributed
  online optimization over unbalanced digraphs via subgradient rescaling,''
  \emph{IEEE Transactions on Control of Network Systems}, 2020.

\bibitem{ding2021differentially}
T.~Ding, S.~Zhu, J.~He, C.~Chen, and X.-P. Guan, ``Differentially private
  distributed optimization via state and direction perturbation in multi-agent
  systems,'' \emph{IEEE Transactions on Automatic Control}, 2021.

\bibitem{nozari2016differentially}
E.~Nozari, P.~Tallapragada, and J.~Cort{\'e}s, ``Differentially private
  distributed convex optimization via functional perturbation,'' \emph{IEEE
  Transactions on Control of Network Systems}, vol.~5, no.~1, pp. 395--408,
  2016.

\bibitem{huang2012differentially}
Z.~Huang, S.~Mitra, and G.~Dullerud, ``Differentially private iterative
  synchronous consensus,'' in \emph{Proceedings of the 2012 ACM workshop on
  Privacy in the electronic society}.\hskip 1em plus 0.5em minus 0.4em\relax
  ACM, 2012, pp. 81--90.

\bibitem{nozari2015}
E.~Nozari, P.~Tallapragada, and J.~Cort\'{e}s., ``Differentially private
  average consensus: obstructions, trade-offs, and optimal algorithm design,''
  \emph{Automatica}, vol.~81, no.~7, pp. 221--231, 2017.

\bibitem{wang2022differentially_consensus}
Y.~Wang, ``A robust dynamic average consensus algorithm that ensures both
  differential privacy and accurate convergence,'' in \emph{62nd IEEE
  Conference on Decision and Control (CDC)}.\hskip 1em plus 0.5em minus
  0.4em\relax IEEE, 2023, pp. 1130--1137.

\bibitem{ye2021differentially}
M.~Ye, G.~Hu, L.~Xie, and S.~Xu, ``Differentially private distributed {Nash}
  equilibrium seeking for aggregative games,'' \emph{IEEE Transactions on
  Automatic Control}, vol.~67, no.~5, pp. 2451--2458, 2021.

\bibitem{wang2022differentially}
J.~Wang, J.-F. Zhang, and X.~He, ``Differentially private distributed
  algorithms for stochastic aggregative games,'' \emph{Automatica}, vol. 142,
  p. 110440, 2022.

\bibitem{camisa2019distributed}
A.~Camisa, F.~Farina, I.~Notarnicola, and G.~Notarstefano, ``Distributed
  constraint-coupled optimization over random time-varying graphs via primal
  decomposition and block subgradient approaches,'' in \emph{2019 IEEE 58th
  Conference on Decision and Control (CDC)}.\hskip 1em plus 0.5em minus
  0.4em\relax IEEE, 2019, pp. 6374--6379.

\bibitem{necoara2013rate}
I.~Necoara and V.~Nedelcu, ``Rate analysis of inexact dual first-order methods
  application to dual decomposition,'' \emph{IEEE Transactions on Automatic
  Control}, vol.~59, no.~5, pp. 1232--1243, 2013.

\bibitem{chang2013distributed}
T.-H. Chang, A.~Nedi{\'c}, and A.~Scaglione, ``Distributed sparse regression by
  consensus-based primal-dual perturbation optimization,'' in \emph{2013 IEEE
  Global Conference on Signal and Information Processing}.\hskip 1em plus 0.5em
  minus 0.4em\relax IEEE, 2013, pp. 289--292.

\bibitem{nedic2009distributed}
A.~Nedi\'c and A.~Ozdaglar, ``Distributed subgradient methods for multi-agent
  optimization,'' \emph{IEEE Transactions on Automatic Control}, vol.~54,
  no.~1, pp. 48--61, 2009.

\bibitem{sundhar2012new}
S.~Sundhar~Ram, A.~Nedi{\'c}, and V.~V. Veeravalli, ``A new class of
  distributed optimization algorithms: Application to regression of distributed
  data,'' \emph{Optimization Methods and Software}, vol.~27, no.~1, pp. 71--88,
  2012.

\bibitem{boyd2004convex}
S.~Boyd and L.~Vandenberghe, \emph{Convex optimization}.\hskip 1em plus 0.5em
  minus 0.4em\relax Cambridge university press, 2004.

\bibitem{bno2003}
D.~P. Berstekas, A.~Nedi\'c, and A.~E. Ozdaglar, \emph{Convex Analysis and
  Optimization}.\hskip 1em plus 0.5em minus 0.4em\relax Athena Scientific,
  2003.

\bibitem{uzawa1958iterative}
H.~Uzawa, ``Iterative methods for concave programming,'' \emph{Studies in
  Linear and Nonlinear Programming}, vol.~6, pp. 154--165, 1958.

\bibitem{kallio1994perturbation}
M.~Kallio and A.~Ruszczynski, ``Perturbation methods for saddle point
  computation,'' \emph{IIASA WP-94-038}, 1994.

\bibitem{kallio1999large}
M.~Kallio and C.~H. Rosa, ``Large-scale convex optimization via saddle point
  computation,'' \emph{Operations Research}, vol.~47, no.~1, pp. 93--101, 1999.

\bibitem{polyak87}
B.~Polyak, ``Introduction to optimization,'' \emph{Optimization Software Inc.,
  Publications Division, New York}, vol.~1, 1987.

\bibitem{dwork2010differential}
C.~Dwork, M.~Naor, T.~Pitassi, and G.~N. Rothblum, ``Differential privacy under
  continual observation,'' in \emph{Proceedings of the forty-second ACM
  Symposium on Theory of Computing}, 2010, pp. 715--724.

\bibitem{chen2023differentially}
X.~Chen, L.~Huang, L.~He, S.~Dey, and L.~Shi, ``A differentially private method
  for distributed optimization in directed networks via state decomposition,''
  \emph{IEEE Transactions on Control of Network Systems}, 2023.

\bibitem{huang2024differential}
L.~Huang, J.~Wu, D.~Shi, S.~Dey, and L.~Shi, ``Differential privacy in
  distributed optimization with gradient tracking,'' \emph{IEEE Transactions on
  Automatic Control}, 2024.

\bibitem{nedic2010constrained}
A.~Nedi\'c, A.~Ozdaglar, and P.~A. Parrilo, ``Constrained consensus and
  optimization in multi-agent networks,'' \emph{IEEE Transactions on Automatic
  Control}, vol.~55, no.~4, pp. 922--938, 2010.

\bibitem{chung1954stochastic}
K.~L. Chung, ``On a stochastic approximation method,'' \emph{The Annals of
  Mathematical Statistics}, pp. 463--483, 1954.

\bibitem{belgioioso2020distributed}
G.~Belgioioso, A.~Nedi{\'c}, and S.~Grammatico, ``Distributed generalized
  {Nash} equilibrium seeking in aggregative games on time-varying networks,''
  \emph{IEEE Transactions on Automatic Control}, vol.~66, no.~5, pp.
  2061--2075, 2021.

\bibitem{nedic2009approximate}
A.~Nedi{\'c} and A.~Ozdaglar, ``Approximate primal solutions and rate analysis
  for dual subgradient methods,'' \emph{SIAM Journal on Optimization}, vol.~19,
  no.~4, pp. 1757--1780, 2009.

\bibitem{bertsekas2015parallel}
D.~Bertsekas and J.~Tsitsiklis, \emph{Parallel and Distributed Computation:
  Numerical Methods}.\hskip 1em plus 0.5em minus 0.4em\relax Athena Scientific,
  2015.

\bibitem{facchinei2003finite}
F.~Facchinei and J.-S. Pang, \emph{Finite-dimensional Variational Inequalities
  and Complementarity Problems}.\hskip 1em plus 0.5em minus 0.4em\relax
  Springer, 2003.

\bibitem{cai2010decentralized}
N.~Cai and J.~Mitra, ``A decentralized control architecture for a microgrid
  with power electronic interfaces,'' in \emph{North American Power Symposium
  2010}.\hskip 1em plus 0.5em minus 0.4em\relax IEEE, 2010, pp. 1--8.

\bibitem{melis2019exploiting}
L.~Melis, C.~Song, E.~De~Cristofaro, and V.~Shmatikov, ``Exploiting unintended
  feature leakage in collaborative learning,'' in \emph{2019 IEEE Symposium on
  Security and Privacy (SP)}.\hskip 1em plus 0.5em minus 0.4em\relax IEEE,
  2019, pp. 691--706.

\bibitem{zhu2019deep}
L.~Zhu, Z.~Liu, and S.~Han, ``Deep leakage from gradients,'' in \emph{Advances
  in Neural Information Processing Systems}, 2019, pp. 14\,774--14\,784.

\end{thebibliography}

\bibliographystyle{IEEEtran}


\end{document}